\documentclass[amsbook,12pt,a4paper]{article}
\usepackage{epsfig, amsfonts, amssymb, amsmath}
\usepackage[usenames]{color}
\usepackage{picinpar}
\usepackage{pgf,pgfarrows,pgfnodes,pgfautomata,pgfheaps,pgfshade}
\newcommand{\nl}{\mbox{}\\}

\begin{document}
\setcounter{page}{1}
%
%
\mbox{} \vspace{-2.000cm} \\
\begin{center}
{\Large \bf
On the supnorm form of Leray's problem} \\
\mbox{} \vspace{-0.300cm} \\
{\Large \bf
for the incompressible Navier-Stokes equations} \\
%
%
%
\nl
\mbox{} \vspace{-0.300cm} \\
{\sc Lineia Sch\"utz, Jana\'\i na P. Zingano and Paulo R. Zingano} \\
\mbox{} \vspace{-0.350cm} \\
{\small Departamento de Matem\'atica Pura e Aplicada} \\
\mbox{} \vspace{-0.670cm} \\
{\small Universidade Federal do Rio Grande do Sul} \\
\mbox{} \vspace{-0.670cm} \\
{\small Porto Alegre, RS 91509, Brazil} \\
\nl
\mbox{} \vspace{-0.300cm} \\
%
%
%
%
{\bf Abstract} \\
\mbox{} \vspace{-0.300cm} \\
\begin{minipage}{12.50cm}
{\small
\mbox{} \hspace{+0.250cm}
We show that
$ {\displaystyle
\;\!
t^{\;\!3/4} \;\!
\|\, \mbox{\boldmath $u$}(\cdot,t) \,
\|_{\scriptstyle L^{\infty}(\mathbb{R}^{3})}
\!\;\!\rightarrow \;\!0
\;\!
} $
as
$ \;\! t \rightarrow \infty \;\!$
for all
Leray-Hopf's global
weak solutions
$ \mbox{\boldmath $u$}(\cdot,t) $
of the incompressible
Navier-Stokes equations \linebreak
in $ \mathbb{R}^{3} \!\:\!$.
It is also shown
that
$ {\displaystyle
\;\!
t \;
\|\, \mbox{\boldmath $u$}(\cdot,t) \;\!-\;\!
e^{\Delta t} \mbox{\boldmath $u$}_0 \,
\|_{\scriptstyle L^{\infty}(\mathbb{R}^{3})}
\!\;\!\rightarrow 0
\;\!
} $
as
$ \;\! t \rightarrow \infty $,
where
$ e^{\Delta t} $
is the heat semigroup,
%
%
as well as other fundamental new results. \linebreak
In spite of the complexity of the questions,
our approach is elementary
and is based on standard tools
like conventional Fourier and energy methods. \\
}
\end{minipage}
\end{center}

\nl
\mbox{} \vspace{-0.550cm} 
\nl
\mbox{} \vspace{-0.800cm} 
\mbox{} \hfill
large time behavior,
Leray's $L^{2}\!\;\!$ problem,
energy estimates, heat semigroup. \\
%
%

%
%

%
%
\nl
\mbox{} \vspace{-0.650cm} \\

{\bf 1. Introduction} \\
\setcounter{section}{1}

In this work,
we derive some new fundamental
large time asymptotic properties
of (globally defined) Leray-Hopf's
weak solutions
\cite{Galdi2000, Leray1934}
of the incompressible Navier-Stokes equations
in three-dimensional space, \\
\mbox{} \vspace{-0.600cm} \\
\begin{equation}
\tag{1.1$a$}
\mbox{\boldmath $u$}_{t} +\,
\mbox{\boldmath $u$} \!\;\!\cdot\!\;\! \nabla\;\!
\mbox{\boldmath $u$}
\,+\,
\nabla p
\,=\;
\Delta \mbox{\boldmath $u$},
\qquad
\nabla \!\cdot\!\;\! \mbox{\boldmath $u$}(\cdot,t)
\,=\,0,
\end{equation}
\mbox{} \vspace{-0.900cm} \\
\begin{equation}
\tag{1.1$b$}
\mbox{\boldmath $u$}(\cdot,0) \,=\,
\mbox{\boldmath $u$}_0 \in L^{2}_{\sigma}(\mathbb{R}^{3}),
\end{equation}
\mbox{} \vspace{-0.200cm} \\
where
$ {\displaystyle
L^{2}_{\sigma}(\mathbb{R}^{3})
} $
denotes the space
of functions
$ {\displaystyle
\;\!
\mbox{\bf u} =
(\:\! \mbox{u}_{\mbox{}_{1}} \!\:\!,
 \:\! \mbox{u}_{\mbox{}_{1}} \!\:\!,
 \:\! \mbox{u}_{\mbox{}_{3}} \!\;\!)
\in L^{2}(\mathbb{R}^{3})
\;\!
} $
with
$ {\displaystyle
\;\!
\nabla \!\cdot \mbox{\bf u}
= 0
} $
in distributional sense.
In his seminal 1934 paper,
Leray \cite{Leray1934}
showed the existence
of (possibly infinitely many) global
weak solutions
$ {\displaystyle
\mbox{\boldmath $u$}(\cdot,t)
\in L^{2}_{\sigma}(\mathbb{R}^{3})
} $
which are
weakly continuous in $L^{2}(\mathbb{R}^{3}) $
and satisfy
$ {\displaystyle
\;\!
\mbox{\boldmath $u$}(\cdot,t)
\in
L^{\infty}([\;\!0, \infty \:\![\:\!, L^{2}(\mathbb{R}^{3}))
\cap
L^{2}([\;\!0, \infty \:\![\:\!, \dot{H}^{1}(\mathbb{R}^{3}))
} $,
with
$ {\displaystyle
\;\!
\|\, \mbox{\boldmath $u$}(\cdot,t) - \mbox{\boldmath $u$}_{0} \;\!
\|_{\mbox{}_{\scriptstyle L^{2}(\mathbb{R}^{3})}}
\!\!\,\!\rightarrow 0
\;\!
} $
as
$ \;\!t \,\mbox{\footnotesize $\searrow$}\;\! 0 \;\!$
and
such that
the energy inequality \\
\mbox{} \vspace{-0.600cm} \\
\begin{equation}
\tag{1.2}
\|\,\mbox{\boldmath $u$}(\cdot,t) \,
\|_{\mbox{}_{\scriptstyle L^{2}(\mathbb{R}^{3})}}^{\:\!2}
\:\!+\;
2 \!
\int_{0}^{\;\!\mbox{\mbox{\footnotesize $t$}}} \!
\|\, D \:\!\mbox{\boldmath $u$}(\cdot,s) \,
\|_{\mbox{}_{\scriptstyle L^{2}(\mathbb{R}^{3})}}^{\:\!2}
\;\!ds
\;\leq\;
\|\,\mbox{\boldmath $u$}_{0} \;\!
\|_{\mbox{}_{\scriptstyle L^{2}(\mathbb{R}^{3})}}^{\:\!2}
\end{equation}
\mbox{} \vspace{-0.100cm} \\
holds
for all $ t \geq 0 $.\footnote{%
%
%
For the definition of the vector norms
involved in (1.2) and other similar expressions
throughout the text,
see (1.16), (1.17) below.
}
%
Moreover,
Leray \cite{Leray1934}
also showed
in his construction
that
there always exists some
$ {\displaystyle
\;\!
t_{\ast} \!\;\!\gg 1
} $
(depending on the solution $\mbox{\boldmath $u$}$)
such that
one actually has
$ {\displaystyle
\mbox{\boldmath $u$} \in
C^{\infty}(\mathbb{R}^{3} \!\times \!\;\![\;\!t_{\ast}, \infty \:\![\:\!)
} $,
$\:\!$and,
for each $ \:\!m \geq 1 $: \\
\mbox{} \vspace{-0.575cm} \\
\begin{equation}
\tag{1.3}
\mbox{\boldmath $u$}(\cdot,t)
\!\;\!\in\!\;\!
L^{\infty}(\:\![\,t_{\ast}, \mbox{\small $T$}\:\!], H^{m}(\mathbb{R}^{3})),
\end{equation}
\mbox{} \vspace{-0.200cm} \\
for each
$ t_{\ast} \!< \mbox{\small $T$} \!< \infty $,
that is,
$ {\displaystyle
\mbox{\boldmath $u$}(\cdot,t)
\!\;\!\in\!\;\!
L^{\infty}_{\tt loc}(\:\![\,t_{\ast}, \infty\:\![\:\!, H^{m}(\mathbb{R}^{3}))
} $.
While the {\em uniqueness} \linebreak
of Leray's solutions remains
a fundamental open question
to this day,
it has been shown by Kato \cite{Kato1984}
and Masuda~\cite{Masuda1984}
(and later by other authors also,
see e.g.$\;$\cite{KajikiyaMiyakawa1986, Wiegner1987})
that {\em all\/} Leray's solutions,
whether uniquely defined by their initial values or not,
must satisfy the
important asymptotic property \\
\mbox{} \vspace{-0.650cm} \\
\begin{equation}
\tag{1.4}
\lim_{t\,\rightarrow\,\infty}
\;\!
\|\, \mbox{\boldmath $u$}(\cdot,t) \,
\|_{\mbox{}_{\scriptstyle L^{2}(\mathbb{R}^{3})}}
=\; 0,
\end{equation}
\mbox{} \vspace{-0.175cm} \\
a question left open in \cite{Leray1934}.
It will prove convenient for our present purposes
that we also provide here a new derivation of (1.4)
along the lines
of the method introduced by
Kreiss, Hagstrom, Lorenz and
one of the authors in
\cite{KreissHagstromLorenzZingano2002, %
KreissHagstromLorenzZingano2003}\footnote{%
%
%
For a detailed account of this method
(mostly due to
T. Hagstrom and J. Lorenz),
see \cite{Moura2013, Rigelo2007}.
}
%
%
to give a straightforward
derivation
of the fundamental
Schonbek-Wiegner decay estimates
\cite{SchonbekWiegner1996, Wiegner1987}
for solutions (and their derivatives)
of the Navier-Stokes equations
in dimension $ n \leq 3 $,
under {\em stronger\/}
assumptions on the initial data.
(See also \cite{OliverTiti2000}.)
It will then be seen
that,
with some extra steps,
one can similarly obtain
the new supnorm result \\
\mbox{} \vspace{-0.075cm} \\
\mbox{} \hspace{+4.200cm}
\fbox{%
\begin{minipage}{5.800cm}
\nl
\mbox{} \vspace{-0.400cm} \\
$ {\displaystyle
\mbox{} \;\;
\lim_{t\,\rightarrow\,\infty}
\;\!
t^{\:\!3/4} \;\!
\|\, \mbox{\boldmath $u$}(\cdot,t) \,
\|_{\mbox{}_{\scriptstyle L^{\infty}(\mathbb{R}^{3})}}
=\; 0
} $, \\
\end{minipage}
}
\mbox{} \vspace{-1.050cm} \\
\mbox{} \hfill (1.5) \\
\nl
\mbox{} \vspace{-0.225cm} \\
%
%
%
which, again, is valid
for all Leray-Hopf's solutions of (1.1),
assuming
$ \mbox{\boldmath $u$}_0 \!\;\!\in L^{2}_{\sigma}(\mathbb{R}^{3}) $
only.
Thus,
by interpolation,
we have,
for any such solution, \\
\mbox{} \vspace{-0.625cm} \\
\begin{equation}
\tag{1.6}
\lim_{t\,\rightarrow\,\infty}
\;\!
t^{\scriptstyle \,\frac{\scriptstyle 3}{\scriptstyle 4}
\,-\, \frac{\scriptstyle 3}{\scriptstyle 2 \:\!q} }
\;\!
\|\, \mbox{\boldmath $u$}(\cdot,t) \,
\|_{\mbox{}_{\scriptstyle L^{q}(\mathbb{R}^{3})}}
=\; 0,
\qquad
2 \leq q \leq \infty,
\end{equation}
\mbox{} \vspace{-0.175cm} \\
uniformly in $q$.
The properties (1.4)$\;\!-\;\!$(1.6)
are well known
and easy to obtain
(see e.g.$\:$\cite{BrazSchutzZingano2013}, Theorem 3.3, p.$\:$95)
for
solutions
$ {\displaystyle
\;\!
\mbox{\boldmath $v$}(\cdot,t) \in
L^{\infty}(\:\![\,t_0, \infty \:\![\:\!, L^{2}(\mathbb{R}^{3}) \:\!)
} $
of the associ\-ated
linear
{\em heat flow\/} problems \\
\mbox{} \vspace{-0.775cm} \\
\begin{equation}
\tag{1.7$a$}
\mbox{\boldmath $v$}_t
\:\!=\:
\Delta
\mbox{\boldmath $v$},
\qquad
\mbox{} \, t > t_0,
\end{equation}
\mbox{} \vspace{-1.050cm} \\
\begin{equation}
\tag{1.7$b$}
\mbox{\boldmath $v$}(\cdot,t_0)
\;\!=\,
\mbox{\boldmath $u$}(\cdot,t_0),
\end{equation}
\mbox{} \vspace{-0.275cm} \\
given
$ \:\!t_0 \!\;\!\geq 0 \:\!$
(arbitrary).
The solution of (1.7)
is given by
$ {\displaystyle
\;\!
\mbox{\boldmath $v$}(\cdot,t)
\:\!=\;\!
e^{\mbox{\scriptsize $\Delta \mbox{\footnotesize $(t - t_0)$}$}}
\mbox{\boldmath $u$}(\cdot,t_0)
} $,
where
$ {\displaystyle
e^{\mbox{\scriptsize $\:\!\Delta \mbox{\footnotesize $\tau$}$}}
\!\!\:\!
} $,
$ \tau \geq 0 $,
denotes the heat semigroup.
It is therefore
natural to think
that
the Leray-Hopf's solutions
of (1.1)
be closely related to
the corresponding heat flows
defined in (1.7).
In fact,
Kato \cite{Kato1984}
obtained
$ {\displaystyle
\lim_{t\,\rightarrow\,\infty}
t^{\:\!1/4 \;-\; \mbox{\footnotesize $\epsilon$}}
\;\!
\|\, \mbox{\boldmath $u$}(\cdot,t) - \mbox{\boldmath $v$}(\cdot,t) \,
\|_{\mbox{}_{\scriptstyle L^{2}(\mathbb{R}^{3})}}
} $
$ = \;\!0 \;\!$
for each
$ \:\!\epsilon > 0 $,
and a bit later
Wiegner \cite{Wiegner1987}
got,
using a very involved argument,
the sharper result\footnote{%
%
%
In addition,
Wiegner obtains (1.4), (1.8)
in the presence of external forces
$ \mbox{\boldmath $f$}(\cdot,t) $,
under suitable assumptions on
$ \mbox{\boldmath $f$}(\cdot,t) $.
Also,
he considers the case
of {\em arbitrary\/} space dimension
$ n \geq 2 $,
which is a complicating factor in the analysis.
While we can certainly extend our approach
to include external forces
$ \mbox{\boldmath $f$} $
in (1.1),
under appropriate assumptions on
$ \mbox{\boldmath $f$}$
which are slightly different from Wiegner's, \linebreak
or similarly extend the analysis down to $ \;\!n = 2 $,
our method is
(as that of
\cite{KreissHagstromLorenzZingano2002, %
KreissHagstromLorenzZingano2003})
limited
to $\;\! n \leq 3 $.
%
}
\mbox{} \\
\mbox{} \vspace{-0.975cm} \\
\begin{equation}
\tag{1.8}
\lim_{t\,\rightarrow\,\infty}
\;\!
t^{\:\!1/4}
\;\!
\|\, \mbox{\boldmath $u$}(\cdot,t) - \mbox{\boldmath $v$}(\cdot,t) \,
\|_{\mbox{}_{\scriptstyle L^{2}(\mathbb{R}^{3})}}
=\; 0
\end{equation}
\mbox{} \vspace{-0.150cm} \\
(see \cite{Wiegner1987}, Theorem ($c$), p.$\;$305).
Again, a simple proof of (1.8)
in the spirit of
\cite{KreissHagstromLorenzZingano2002, %
KreissHagstromLorenzZingano2003} \linebreak
is pro\-vided here (see Section~3),
after (1.4) has been obtained.
This is useful
to pave our way
for the corresponding
sup\-norm result
obtained in Section~4,
viz., \\
\mbox{} \vspace{-0.075cm} \\
\mbox{} \hspace{+3.670cm}
\fbox{%
\begin{minipage}{6.850cm}
\nl
\mbox{} \vspace{-0.400cm} \\
$ {\displaystyle
\mbox{} \;\;
\lim_{t\,\rightarrow\,\infty}
\;\!
t
\:
\|\, \mbox{\boldmath $u$}(\cdot,t) - \mbox{\boldmath $v$}(\cdot,t) \,
\|_{\mbox{}_{\scriptstyle L^{\infty}(\mathbb{R}^{3})}}
=\; 0
} $. \\
\end{minipage}
}
\mbox{} \vspace{-1.050cm} \\
\mbox{} \hfill (1.9) \\
\nl
\mbox{} \vspace{-0.200cm} \\
%
%
%
%
By interpolation,
it follows
from (1.8), (1.9) that \\
\mbox{} \vspace{-0.625cm} \\
\begin{equation}
\tag{1.10}
\lim_{t\,\rightarrow\,\infty}
\;\!
t^{\:\!1 \,-\, \frac{\scriptstyle 3}{\scriptstyle \:\!2 \:\!q\:\!}}
\;\!
\|\, \mbox{\boldmath $u$}(\cdot,t) \,-\,
     \mbox{\boldmath $v$}(\cdot,t) \,
\|_{\mbox{}_{\scriptstyle L^{q}(\mathbb{R}^{3})}}
=\; 0,
\qquad
2 \leq q \leq \infty,
\end{equation}
\mbox{} \vspace{-0.200cm} \\
uniformly in $q$.
It is worth noticing
that
these results
improve the
previous estimates \\
\mbox{} \vspace{-0.600cm} \\
\begin{equation}
\tag{1.11$a$}
\limsup_{t\,\rightarrow\,\infty}
\;
t^{\scriptstyle \,\frac{\scriptstyle 3}{\scriptstyle 4}
\,-\, \frac{\scriptstyle 3}{\scriptstyle 2 \:\!q} }
\;\!
\|\, \mbox{\boldmath $u$}(\cdot,t) \,
\|_{\mbox{}_{\scriptstyle L^{q}(\mathbb{R}^{3})}}
<\, \infty,
\end{equation}
%
%
\begin{equation}
\tag{1.11$b$}
\limsup_{t\,\rightarrow\,\infty}
\;
t^{\:\!1 \,-\, \frac{\scriptstyle 3}{\scriptstyle \:\!2 \:\!q\:\!}}
\;\!
\|\, \mbox{\boldmath $u$}(\cdot,t) \,-\,
     \mbox{\boldmath $v$}(\cdot,t) \,
\|_{\mbox{}_{\scriptstyle L^{q}(\mathbb{R}^{3})}}
<\, \infty
\end{equation}
\mbox{} \vspace{-0.150cm} \\
obtained
by Beir\~ao da Veiga and Wiegner
in \cite{BeiraoDaVeiga1987, Wiegner1990}
for finite $ \:\! q > 2 $.
Here is a brief overview of what is next.
After some important mathematical preliminaries
on the Leray-Hopf's solutions
to the Navier-Stokes system (1.1)
have been reviewed in Section 2
for later use,
along with two new fundamental results
given by Theorems 2.2 and 2.3,
we turn our attention to the basic
$L^{2}\!\;\!$ estimates (1.4) and (1.8),
which are rederived in Section 3
along the lines of
\cite{KreissHagstromLorenzZingano2002, KreissHagstromLorenzZingano2003}.
This shows the way
to obtain the
more difficult estimates
(1.5) and (1.9),
which is the goal of Section 4.
In these two sections,
the key point is to
first observe
that \\
\mbox{} \vspace{-0.640cm} \\
\begin{equation}
\tag{1.12}
\lim_{t\,\rightarrow\,\infty}
\;\!
t^{\:\!1/2} \;\!
\|\, D \mbox{\boldmath $u$}(\cdot,t) \,
\|_{\mbox{}_{\scriptstyle L^{2}(\mathbb{R}^{3})}}
\!\;\!=\; 0,
\end{equation}
\mbox{} \vspace{-0.190cm} \\
from which the desired estimates
can be more easily obtained.
%
%
Although we restrict our attention here
to dimension $ n = 3 $,
it will be clearly seen that
the method can also be used
in the case $ n = 2 $,
which is actually easier
since (1.12) turns out to be trivial
in this case.
Put together,
the results
for $ \;\!n = 2, 3 \;\!$
can be summarized
as follows.
One has,
for each
$ \:\!2 \leq q \leq \infty $
(and $ \:\! n = 2, 3 $): \\
\mbox{} \vspace{-0.600cm} \\
\begin{equation}
\tag{1.13$a$}
\lim_{t\,\rightarrow\,\infty}
\;\!
t^{\scriptstyle \;\!\frac{\scriptstyle n}{\scriptstyle 4}
\,-\, \frac{\scriptstyle n}{\scriptstyle 2 \:\!q} }
\;\!
\|\, \mbox{\boldmath $u$}(\cdot,t) \,
\|_{\mbox{}_{\scriptstyle L^{q}(\mathbb{R}^{n})}}
=\; 0,
\end{equation}
\mbox{} \vspace{-0.770cm} \\
\begin{equation}
\tag{1.13$b$}
\lim_{t\,\rightarrow\,\infty}
\;\!
t^{\;\! \frac{\scriptstyle n - 1}{\scriptstyle 2}
\,-\, \frac{\scriptstyle n}{\scriptstyle \:\!2 \:\!q\:\!}}
\;\!
\|\, \mbox{\boldmath $u$}(\cdot,t) \,-\,
     \mbox{\boldmath $v$}(\cdot,t) \,
\|_{\mbox{}_{\scriptstyle L^{q}(\mathbb{R}^{n})}}
=\; 0,
\end{equation}
\mbox{} \vspace{-0.250cm} \\
uniformly in $\:\!q \in [\:\!2, \infty\:\!] $,
where
$ {\displaystyle
\;\!
\mbox{\boldmath $v$}(\cdot,t)
=\:\!
e^{\mbox{\scriptsize $\Delta \mbox{\footnotesize $(t - t_0)$}$}}
\mbox{\boldmath $u$}(\cdot,t_0)
} $,
$ \:\! t_0 \!\geq 0 $ arbitrary,
see (1.7),
under the sole assumption
that
$ \mbox{\boldmath $u$}_0 \!\in L^{2}_{\sigma}(\mathbb{R}^{n}) $.
A proof (or disproof)
of this general property
in higher dimensions
is apparently still missing
in the literature.
For $ n \leq 3 $, \linebreak
everything needed to obtain (1.13)
was already known by 1934
after the publication
of \cite{Leray1934},
as the next sections show
--- and yet it has taken full fifty years
before even the easier part of (1.13)
could have finally been established!
We hope that this shows the power
of the ideas presented here,
as well as of the approach
introduced in
\cite{KreissHagstromLorenzZingano2002, %
KreissHagstromLorenzZingano2003}. \linebreak
In fact,
a deeper combination
of these ideas
has now led
to the complete solution
of the full Leray's problem
in dimension $ n \leq 3 $
\cite{HagstromLorenzZingano2015}:
one has,
for every
$ s \geq 0 $,
and
any
$ 0 \leq t_0 \leq t_1 \!\;\!< t $, \\
\mbox{} \vspace{-0.900cm} \\
\begin{equation}
\tag{1.14$a$}
\lim_{t\,\rightarrow\,\infty}
\,
t^{\mbox{}^{\scriptstyle \frac{\scriptstyle s}{2} }}
\|\, \mbox{\boldmath $u$}(\cdot,t) \,
\|_{\mbox{}_{\scriptstyle \dot{H}^{\!\:\!s}\!\;\!(\mathbb{R}^{n})}}
\!\;\!=\:0,
\end{equation}
\mbox{} \vspace{-0.710cm} \\
\begin{equation}
\tag{1.14$b$}
\lim_{t\,\rightarrow\,\infty}
\,
t^{\mbox{}^{\scriptstyle
\frac{\scriptstyle n}{4} \;\!+\;\!
\frac{\scriptstyle s}{2} \;\!-\;\!\frac{1}{2}
}}
\;\!
\|\, \mbox{\boldmath $u$}(\cdot,t) \:\!-\:\!
e^{\:\!\Delta (t - t_0)} \mbox{\boldmath $u$}(\cdot,t_0) \,
\|_{\mbox{}_{\scriptstyle \dot{H}^{\!\:\!s}\!\;\!(\mathbb{R}^{n})}}
\!\;\!=\:0,
\end{equation}
\mbox{} \vspace{-0.850cm} \\
\begin{equation}
\tag{1.14$c$}
\|\, e^{\:\!\Delta (t - t_0)}
\mbox{\boldmath $u$}(\cdot,t_0) \;\!-\;\!
e^{\:\!\Delta (\:\!t - t_1\!\;\!)}
\mbox{\boldmath $u$}(\cdot,t_1) \,
\|_{\mbox{}_{\scriptstyle \dot{H}^{\!\:\!s}\!\;\!(\mathbb{R}^{n})}}
\!\:\!\leq\:\!
K\!\:\!(n,s) \,
(\:\! t_1 \!\;\!-\;\! t_0 )^{\mbox{}^{\scriptstyle \!\frac{1}{2}}}
(\:\! t \:\!-\;\! t_{1} \!\;\!)^{\mbox{}^{\scriptstyle \!\!\!\;\!-\,
\left( \frac{\scriptstyle n}{4} \;\!+\;\! \frac{\scriptstyle s}{2}
\right)}}
\end{equation}
\mbox{} \vspace{-0.150cm} \\
for arbitrary Leray solutions
in $ \mathbb{R}^{n}\!$,
under the unique assumption
of square-integrable, divergence-free
initial data.
(For the more involved analysis
giving (1.14),
see~\cite{HagstromLorenzZingano2015}.
In the simpler case of
dimension $ n = 2 $,
(1.14$a$) was shown in
\cite{BenameurSelmi2012}
by a different method. \linebreak
Some related high-order estimates
have also been obtained in
\cite{KreissHagstromLorenzZingano2002, %
OliverTiti2000, SchonbekWiegner1996},
but under \linebreak
stronger assumptions on the initial data.)
$\!\:\!$Here,
$ \!\:\!\dot{H}^{\!\;\!s}\!\;\!(\mathbb{R}^{n}) $
denotes the homogeneous Sobolev space
of all functions
$ {\displaystyle
\mbox{\boldmath $v$} = (v_{\scriptscriptstyle 1}\!\;\!,\!...,v_{n})
\in L^{2}(\mathbb{R}^{n})
} $
such that
$ {\displaystyle
\;\!
| \cdot |^{\,\!s} \:\!
|\,\hat{\mbox{\boldmath $v$}}(\cdot)\,|
\in\!\;\! L^{2}(\mathbb{R}^{n})
} $,
where
$ \hat{\mbox{\boldmath $v$}}(\cdot) $
stands for the Fourier transform of
$ \mbox{\boldmath $v$}(\cdot) $,
with norm
$ {\displaystyle
\,\!
\| \cdot
\|_{\mbox{}_{\scriptstyle \dot{H}^{\!\:\!s}\!\;\!(\mathbb{R}^{n})}}
\!
} $
defined by \\
\mbox{} \vspace{-0.525cm} \\
\begin{equation}
\tag{1.15}
\|\, \mbox{\boldmath $v$} \,
\|_{\mbox{}_{\scriptstyle \dot{H}^{\!\:\!s}\!\;\!(\mathbb{R}^{n})}}
=\,
\Bigl\{%
\int_{\mathbb{R}^{n}} \!\!
|\,\xi\,|_{\mbox{}_{2}}^{\:\!2\:\!s}
\;\!
|\, \hat{\mbox{\boldmath $v$}}(\xi)\,|_{\mbox{}_{2}}^{2}
\,d\xi
\:\Bigr\}^{\!\!\;\!1/2}
\!\!.
\end{equation}
\mbox{} \vspace{-0.100cm} \\
Thus,
one has
$ {\displaystyle
\,
t^{\frac{\scriptstyle m}{2} }
\,\!
\|\, D^{m} \mbox{\boldmath $u$}(\cdot,t) \,
\|_{\mbox{}_{\scriptstyle L^{2}(\mathbb{R}^{n})}}
\!\!\;\!\rightarrow 0
} $,
$ {\displaystyle
\;\!
t^{\frac{\scriptstyle n}{4} + \frac{\scriptstyle m}{2} - \frac{1}{2} }
\,\!
\|\, D^{m} \bigl\{
\mbox{\boldmath $u$}(\cdot,t) -\:\!
e^{\:\!\Delta t} \mbox{\boldmath $u$}_0
\bigr\}\:\!
\|_{\mbox{}_{\scriptstyle L^{2}(\mathbb{R}^{n})}}
\!\rightarrow 0
} $
(as $ t \rightarrow \infty $)
for every $ m \in \mathbb{N} $,
and so forth,
which are important extensions of (1.12).
In (1.15) and throughout the text,
$ | \cdot |_{\mbox{}_{2}} $
denotes the Euclidean norm
in $ \mathbb{R}^{n} \!$. \\
\mbox{} \vspace{-0.950cm} \\

More on notation:
boldface letters
are used for
vector quantities,
as in
$ {\displaystyle
\;\!
\mbox{\boldmath $u$}(x,t)
=
} $
$ {\displaystyle
(\:\! u_{\mbox{}_{\!\:\!1}}\!\;\!(x,t),
 \:\! u_{\mbox{}_{\!\;\!2}}\!\;\!(x,t),
 \:\! u_{\mbox{}_{\!\;\!3}}\!\;\!(x,t) \:\!)
} $.
$\!$Also,
$ \nabla p \;\!\equiv \nabla p(\cdot,t) $
denotes the spatial gradient of $ \;\!p(\cdot,t) $,
$ D_{\!\;\!j} \!\;\!=\:\! \partial / \partial x_{\!\;\!j} \!\;\! $,
and
$ {\displaystyle
\:\!
\nabla \!\cdot \mbox{\boldmath $u$}
\:\!=
  D_{\mbox{}_{\!\:\!1}} u_{\mbox{}_{\!\:\!1}} \!\;\!+
  D_{\mbox{}_{\!\;\!2}} \:\! u_{\mbox{}_{\!\;\!2}} \!\;\!+
  D_{\mbox{}_{\!\;\!3}} \:\! u_{\mbox{}_{\!\;\!3}} )
} $
is the (spatial) divergence of
$ \:\!\mbox{\boldmath $u$}(\cdot,t) $.
$ {\displaystyle
\| \;\!\cdot\;\!
\|_{\mbox{}_{\scriptstyle L^{q}(\mathbb{R}^{3})}}
\!\;\!
} $,
$ 1 \leq q \leq \infty $,
denote the standard norms
of the Lebesgue spaces
$ L^{q}(\mathbb{R}^{3}) $,
with \\
\mbox{} \vspace{-1.050cm} \\
\begin{equation}
\tag{1.16$a$}
\|\, \mbox{\boldmath $u$}(\cdot,t) \,
\|_{\mbox{}_{\scriptstyle L^{q}(\mathbb{R}^{3})}}
\;\!=\;
\Bigl\{\,
\sum_{i\,=\,1}^{3} \int_{\mathbb{R}^{3}} \!
|\:u_{i}(x,t)\,|^{q} \;\!dx
\,\Bigr\}^{\!\!\:\!1/q}
\end{equation}
\mbox{} \vspace{-0.750cm} \\
\begin{equation}
\tag{1.16$b$}
\|\, D \mbox{\boldmath $u$}(\cdot,t) \,
\|_{\mbox{}_{\scriptstyle L^{q}(\mathbb{R}^{3})}}
\;\!=\;
\Bigl\{\,
\sum_{i, \,j \,=\,1}^{3} \int_{\mathbb{R}^{3}} \!
|\, D_{\!\;\!j} \;\!u_{i}(x,t)\,|^{q} \;\!dx
\,\Bigr\}^{\!\!\:\!1/q}
\end{equation}
\mbox{} \vspace{-0.750cm} \\
\begin{equation}
\tag{1.16$c$}
\|\, D^{2} \mbox{\boldmath $u$}(\cdot,t) \,
\|_{\mbox{}_{\scriptstyle L^{q}(\mathbb{R}^{3})}}
\;\!=\;
\Bigl\{\!\!
\sum_{\mbox{} \;\;i, \,j, \,\ell \,=\,1}^{3}
\!\;\! \int_{\mathbb{R}^{3}} \!
|\, D_{\!\;\!j} \:\!D_{\ell} \, u_{i}(x,t)\,|^{q} \;\!dx
\,\Bigr\}^{\!\!\:\!1/q}
\end{equation}
\mbox{} \vspace{-0.050cm} \\
if $ 1 \leq q < \infty $,
and
$ {\displaystyle
\;\!
\|\, \mbox{\boldmath $u$}(\cdot,t) \,
\|_{\mbox{}_{\scriptstyle L^{\infty}(\mathbb{R}^{3})}}
\!=\;\!
\max \, \bigl\{\,
\|\,u_{i}(\cdot,t)\,
\|_{\mbox{}_{\scriptstyle L^{\infty}(\mathbb{R}^{3})}}
\!\!: \, 1 \leq i \leq 3
\,\bigr\}
} $
if $ \;\!q = \infty $.
We will also find it convenient
in many places
to use the following alternative definition
for the supnorm of $ \mbox{\boldmath $u$}(\cdot,t) $: \\
\mbox{} \vspace{-0.675cm} \\
\begin{equation}
\tag{1.17}
\|\, \mbox{\boldmath $u$}(\cdot,t) \,
\|_{\mbox{}_{\scriptstyle \infty}}
=\;
\mbox{ess}\,\sup\; \bigl\{\:
|\,\mbox{\boldmath $u$}(x,t) \,|_{\mbox{}_{2}}
\!\!\:\!: \: x \in \mathbb{R}^{3}
\,\bigr\}.
\end{equation}
\mbox{} \vspace{-0.250cm} \\
All other notation, when not standard,
will be explained as it appears in the text. \linebreak
\mbox{} \vspace{-0.950cm} \\
%

For readers interested mainly
in the new results
obtained in the present work,
one could at this point go directly to
Theorems 2.2 and 2.3 in Section 2,
and
Theorems 4.1 and 4.2 in Section 4,
with a quick pass at
(2.22) and the
discussions in Section~3
and the \mbox{\small \sc Appendix},
particularly Theorem A.1.
The few remaining results
may be also
worth browsing,
as some are not so widely known
as they surely deserve to be.
%
%

%
%

%
%
$\mbox{\boldmath $\S\;\!2.$}$
{\bf Some mathematical preliminaries} \\
\setcounter{section}{2}
\mbox{} \vspace{-0.450cm} \\

In this section,
we collect some basic results
that will play an important role
later \linebreak
in our derivation
of (1.4), (1.5), (1.8) and (1.9),
and we also introduce
two fundamental new results
(Theorems 2.2 and 2.3 below).
For the construction of
Leray-Hopf's
solutions
$ \mbox{\boldmath $u$}(\cdot,t) $
to the Navier-Stokes equations (1.3),
see e.g.$\;$\cite{Galdi2000, Leray1934}.
These solutions \linebreak
were originally
obtained in \cite{Leray1934}
by introducing an ingenious
regularization procedure
which, for convenience,
is briefly reviewed next.
Taking
(any)
$ {\displaystyle
\;\!
G \in C^{\infty}_{0}(\mathbb{R}^{n})
\;\!
} $
nonnegative
with
$ \!\;\!\int_{\mathbb{^R}^{3}} \!\:\!G(x) \,dx \:\!=\:\! 1 $
and
setting
$ {\displaystyle
\;\!
\bar{\mbox{\boldmath $u$}}_{\mbox{}_{\scriptstyle \!\;\!0, \,\delta}}
\!\;\!(\cdot)
\in C^{\infty}(\mathbb{R}^{3})
} $
by
convolving
$ {\displaystyle
\;\!
\mbox{\boldmath $u$}_{0}(\cdot)
\;\!
} $
with
$ {\displaystyle
\;\!
G_{\mbox{}_{\scriptstyle \!\delta}}(x)
\:\!=\:\!
\delta^{\;\!-\,n} \;\!G(x/\delta)
} $,
$ \;\!\delta > 0 $,
if
we define
$ {\displaystyle
\;\!
\mbox{\boldmath $u$}_{\mbox{}_{\scriptstyle \!\;\!\delta}},
\,
p_{\mbox{}_{\scriptstyle \!\;\!\delta}}
\in
C^{\infty}(\:\!\mathbb{R}^{3} \!\times\!\;\! [\,0, \infty\:\![\:\!)
} $
as the
(unique, globally defined)
classical
$L^{2}$ solutions
of the associated equations \\
\mbox{} \vspace{-0.575cm} \\
\begin{equation}
\tag{2.1$a$}
\mbox{\small $ {\displaystyle \frac{\partial}{\partial \;\!t} }$}
\,\mbox{\boldmath $u$}_{\mbox{}_{\scriptstyle \!\;\!\delta}}
\:\!+\:
\bar{\mbox{\boldmath $u$}}_{\mbox{}_{\scriptstyle \!\;\!\delta}}
(\cdot,t) \!\;\!
\cdot\!\;\! \nabla\;\!
\mbox{\boldmath $u$}_{\mbox{}_{\scriptstyle \!\;\!\delta}}
\:\!+\,
\nabla \:\!p_{\mbox{}_{\scriptstyle \!\;\!\delta}}
\;=\;
\Delta \:\!\mbox{\boldmath $u$}_{\mbox{}_{\scriptstyle \!\;\!\delta}},
\qquad
\nabla \!\cdot\!\;\!
\mbox{\boldmath $u$}_{\mbox{}_{\scriptstyle \!\;\!\delta}}(\cdot,t)
\,=\,0,
\end{equation}
\mbox{} \vspace{-0.950cm} \\
\begin{equation}
\tag{2.1$b$}
\mbox{\boldmath $u$}_{\mbox{}_{\scriptstyle \!\;\!\delta}}(\cdot,0)
\,=\,
\bar{\mbox{\boldmath $u$}}_{\mbox{}_{\scriptstyle \!\;\!0, \,\delta}}
\!\;\!:=\,
G_{\mbox{}_{\scriptstyle \!\delta}}
\!\:\!\ast
\mbox{\boldmath $u$}_{0}
\;\!
\in\!
\bigcap_{m\,=\,1}^{\infty}
\! H^{m}(\mathbb{R}^{3}),
\end{equation}
\mbox{} \vspace{-0.100cm} \\
where
$ {\displaystyle
\;\!
\bar{\mbox{\boldmath $u$}}_{\mbox{}_{\scriptstyle \!\;\!\delta}}
\!\;\!(\cdot,t)
\!\;\!:=\:\!
G_{\mbox{}_{\scriptstyle \!\delta}} \!\,\!\ast
{\mbox{\boldmath $u$}}_{\mbox{}_{\scriptstyle \!\;\!\delta}}
\!\;\!(\cdot,t)
} $,
it was shown by Leray
that,
for some sequence
$ {\displaystyle
\;\!
\delta^{\;\!\prime}
\!\rightarrow
\:\!0
} $,
one has
the weak convergence property \\
\mbox{} \vspace{-0.600cm} \\
\begin{equation}
\tag{2.2}
\mbox{\boldmath $u$}_{\mbox{}_{\scriptstyle \!\;\!\delta^{\;\!\prime}}}(\cdot,t)
\,\rightharpoonup \,
\mbox{\boldmath $u$}(\cdot,t)
\quad \;\,
\mbox{as } \;\,\delta^{\;\!\prime}
\!\rightarrow
\:\!0,
\qquad \;\;\;
\forall \;\, t \geq 0,
\end{equation}
\mbox{} \vspace{-0.200cm} \\
that is,
$ {\displaystyle
\mbox{\boldmath $u$}_{\mbox{}_{\scriptstyle \!\;\!\delta^{\;\!\prime}}}(\cdot,t)
\,\rightarrow \,
\mbox{\boldmath $u$}(\cdot,t)
\;\!
} $
weakly in $L^{2}(\mathbb{R}^{3}) $,
for every $ \;\! t \geq 0 $
(see \cite{Leray1934}, p.$\;$237),
with \linebreak
$ {\displaystyle
\mbox{\boldmath $u$}(\cdot,t)
\in
L^{\infty}([\;\!0, \infty \:\![\:\!, L^{2}_{\sigma}(\mathbb{R}^{3}))
\cap
L^{2}([\;\!0, \infty \:\![\:\!, \mbox{$\stackrel{.}{H}$}\mbox{}^{1}(\mathbb{R}^{3}))
\cap
C^{0}_{w}([\;\!0, \infty \:\![\:\!, L^{2}(\mathbb{R}^{3}))
\:\!
} $
continuous in $L^{2}\!\;\! $
at $ t = 0 $
and
solving
the Navier-Stokes equations (1.1$a$)
in distributional sense.
Moreover,
the energy inequality (1.2)
is satisfied
for all $ t \geq 0 $,
so that,
in particular, \\
\mbox{} \vspace{-0.475cm} \\
\begin{equation}
\tag{2.3}
\int_{0}^{\:\!\infty} \!\!\!
\|\, D \mbox{\boldmath $u$}(\cdot,t) \,
\|_{\mbox{}_{\scriptstyle L^{2}(\mathbb{R}^{3})}}^{\:\!2}
\:\!dt
\;\leq\;
\mbox{\small $ {\displaystyle
\frac{1}{\:\!2\:\!} }$} \:
\|\: \mbox{\boldmath $u$}_0 \,
\|_{\mbox{}_{\scriptstyle L^{2}(\mathbb{R}^{3})}}^{\:\!2}
\!.
\end{equation}
\mbox{} \vspace{-0.150cm} \\
Another important property
shown in \cite{Leray1934}
is that
$ {\displaystyle
\;\!
\mbox{\boldmath $u$} \in C^{\infty}(\:\![\,t_{\!\;\!\ast}, \infty\:\![\;\!)
\;\!
} $
for some $ \;\! t_{\!\;\!\ast} \!\gg\!\;\! 1 $,
with
$ {\displaystyle
D^{m} \mbox{\boldmath $u$}(\cdot,t)
\in L^{\infty}_{\tt loc} (\:\![\, t_{\!\;\!\ast}, \infty \:\![\:\!, L^{2}(\mathbb{R}^{3})\:\!)
\;\!
} $
for each $ \:\! m \geq 1 $.
This fact
(together with Theorems 2.2 and 2.3 below)
will greatly simplify
our analysis in Sections 3 and 4. \linebreak
Other results
needed later
have mostly to do with the
Helmholtz-Weyl projection
of
$ {\displaystyle
-\,
\mbox{\boldmath $u$}(\cdot,t)
\!\;\!\cdot\!\;\! \nabla \:\!
\mbox{\boldmath $u$}(\cdot,t)
} $
into
$ L^{2}_{\sigma}(\mathbb{R}^{3}) $,
that is,
the divergence-free field
$ \:\!\mbox{\boldmath $Q$}(\cdot,t) \in L^{2}(\mathbb{R}^{3}) \:\!$
given by \\
\mbox{} \vspace{-0.900cm} \\
\begin{equation}
\tag{2.4}
\mbox{\boldmath $Q$}(\cdot,t)
\;\! := \;
-\:
\mbox{\boldmath $u$}(\cdot,t)
\!\;\!\cdot\!\;\! \nabla \:\!
\mbox{\boldmath $u$}(\cdot,t)
\,-\;\! \nabla \:\!p\:\!(\cdot,t),
\qquad
\mbox{a.e. }\;\! t > 0.
\end{equation}
\mbox{} \vspace{-0.250cm} \\
For convenience,
they are discussed in more detail
in the remainder of this section.
\mbox{} \vspace{-0.500cm} \\
%
%
%
%
{\bf Theorem 2.1.}
\textit{%
For almost every $ \;\! s > 0 $,
one has
} \\
\mbox{} \vspace{-0.700cm} \\
\begin{equation}
\tag{2.5$a$}
\|\: e^{\Delta (\:\!\mbox{\footnotesize $t$} \;\!-\, \mbox{\footnotesize $s$})}
\:\! \mbox{\boldmath $Q$}(\cdot,s) \,
\|_{\mbox{}_{\scriptstyle L^{2}(\mathbb{R}^{3})}}
\;\!\leq\:
K \;\!
( t - s )^{-\,3/4} \,
\|\, \mbox{\boldmath $u$}(\cdot,s) \,
\|_{\mbox{}_{\scriptstyle L^{2}(\mathbb{R}^{3})}}
\:\!
\|\, D \mbox{\boldmath $u$}(\cdot,s) \,
\|_{\mbox{}_{\scriptstyle L^{2}(\mathbb{R}^{3})}}
\end{equation}
\mbox{} \vspace{-0.400cm} \\
\textit{and} \\
\mbox{} \vspace{-0.900cm} \\
\begin{equation}
\tag{2.5$b$}
\|\: e^{\Delta (\:\!\mbox{\footnotesize $t$} \;\!-\, \mbox{\footnotesize $s$})}
\:\! \mbox{\boldmath $Q$}(\cdot,s) \,
\|_{\mbox{}_{\scriptstyle \infty}}
\,\leq\:
K \;\!
( t - s )^{-\,3/4} \,
\|\, \mbox{\boldmath $u$}(\cdot,s) \,
\|_{\mbox{}_{\scriptstyle \infty}}
\:\!
\|\, D \mbox{\boldmath $u$}(\cdot,s) \,
\|_{\mbox{}_{\scriptstyle L^{2}(\mathbb{R}^{3})}}
\end{equation}
\mbox{} \vspace{-0.200cm} \\
\textit{%
for all $ \,t > s $,
where
$ \:\!K \!\:\!=\:\! (\:\! 8 \:\!\pi )^{-\,3/4} \!\:\!$.
} \\
%
%
\nl
%
%
{\small
{\bf Proof:}
The following argument
is adapted from
\cite{KreissHagstromLorenzZingano2002}.
Considering (2.5$a$) first,
let
$ {\displaystyle
\;\!
\mathbb{F}\:\![\;\!f\;\!]
\equiv
\hat{f}
\;\!
} $
denote the Fourier transform
of a given function
$ \!\;\!f \!\;\!\in L^{1}(\mathbb{R}^{3}) $,
viz., \\
\mbox{} \vspace{-0.525cm} \\
\begin{equation}
\tag{2.6}
\mathbb{F}\:\![\;\!f\;\!]\:\!(k)
\:\equiv\;
\hat{f}(k)
\,:=\;
(\:\!2\:\!\pi)^{-\,3/2} \!\!
\int_{\mbox{}_{\scriptstyle \!\:\!\mathbb{R}^{3}}}
\!\!\:\!
e^{\mbox{\scriptsize $- \!\;\!
\stackrel{\mbox{\tiny $\circ$}}{\mbox{\sf \i$\!\;\!$\i}}
\!\:\! k \!\;\!\cdot\!\;\! x$}}
f(x) \: dx,
\qquad
k \in \mathbb{R}^{3}
%
%
\end{equation}
\mbox{} \vspace{-0.150cm} \\
(where
$ {\displaystyle
\:\!
\stackrel{\mbox{\tiny $\circ$}}{\mbox{\sf \i$\!\;\!$\i}}
\mbox{}\!\!\;\!\mbox{}^{2}
\!\;\!=\;\! - \;\!1
} $).
Given
$ {\displaystyle
\;\!
\mbox{\bf v}(\cdot,s) =
(\:\! \mbox{v}_{\mbox{}_{\!\;\!1}}\!\;\!(\cdot,s),
      \mbox{v}_{\mbox{}_{\!\;\!2}}\!\;\!(\cdot,s),
      \mbox{v}_{\mbox{}_{\!\;\!3}}\!\;\!(\cdot,s)\:\! )
\in
L^{1}(\mathbb{R}^{3}) \cap L^{2}(\mathbb{R}^{3})
\:\!
} $
arbitrary,
we~get,
using Parseval's identity, \\
\mbox{} \vspace{-0.520cm} \\
\begin{equation}
\notag
\begin{split}
\|\:  e^{\:\!\Delta \mbox{\scriptsize $(t-s)$}} \;\!
\mbox{\bf v}(\cdot,s) \,
\|_{\mbox{}_{\scriptstyle L^{2}(\mathbb{R}^{3})}}^{\:\!2}
\;\!&=\;\;\!
\|\; \mathbb{F}\;\![\, e^{\:\!\Delta \mbox{\scriptsize $(t-s)$}}
\;\!\mbox{\bf v}(\cdot,s) \: ] \,
\|_{\mbox{}_{\scriptstyle L^{2}(\mathbb{R}^{3})}}^{\:\!2} \\
&=\,
\int_{\mathbb{R}^{3}} \!\!
e^{-\,2\,|\,\mbox{\scriptsize $k$}\,
|_{\mbox{}_{\!\;\!2}}^{\:\!\scriptstyle 2} \:\!(t \,-\, s)}
\;\!
|\, \hat{\mbox{\bf v}}(k,s) \,
|_{\mbox{}_{\!\;\!2}}^{\:\!2}
\, dk \\
&\leq\;
\|\, \hat{\mbox{\bf v}}(\cdot,s) \,
\|_{\mbox{}_{\scriptstyle \infty}}^{\:\!2} \!\!\;\!
\int_{\mathbb{R}^{3}} \!\!
e^{-\,2\,|\,\mbox{\scriptsize $k$}\,
|_{\mbox{}_{\!\;\!2}}^{\:\!\scriptstyle 2} \:\!(t \,-\, s)}
\, dk \\
&=\;
\Bigl(\;\! \frac{\;\!\pi\;\!}{2} \,\Bigr)^{\!\!\;\!3/2}
\!\;\!
( t - s )^{- \,3/2}
\,
\|\, \hat{\mbox{\bf v}}(\cdot,s) \,
\|_{\mbox{}_{\scriptstyle \infty}}^{\:\!2}
\!\;\!,
\end{split}
\end{equation}
\mbox{} \vspace{-0.300cm} \\
that is, \\
\mbox{} \vspace{-0.800cm} \\
\begin{equation}
\tag{2.7}
\|\:  e^{\:\!\Delta \mbox{\scriptsize $(t-s)$}} \;\!
\mbox{\bf v}(\cdot,s) \,
\|_{\mbox{}_{\scriptstyle L^{2}(\mathbb{R}^{3})}}
\;\!\leq\;\;\!
\Bigl(\;\! \frac{\;\!\pi\;\!}{2} \,\Bigr)^{\!\!\;\!3/4}
\!\;\!
( t - s )^{- \,3/4}
\,
\|\, \hat{\mbox{\bf v}}(\cdot,s) \,
\|_{\mbox{}_{\scriptstyle \infty}}
\!\;\!,
\end{equation}
\mbox{} \vspace{-0.175cm} \\
where
$ {\displaystyle
\;\!
|\;\!\cdot\;\!|_{\mbox{}_{2}}
} $
denotes the Euclidean norm
in $ \mathbb{R}^{3}\!\;\!$
and
$ {\displaystyle
\;\!
\|\, \hat{\mbox{\bf v}}(\cdot,s) \,
\|_{\mbox{}_{\scriptstyle \infty}}
\!=\,
\sup \;\{\:
|\, \hat{\mbox{\bf v}}(k,s) \,|_{\mbox{}_{2}}
\!: \, k\,\in\,\mathbb{R}^{3} \:\!\}
} $.
As will be shown next,
(2.5$a$) follows from
a direct application of (2.7)
to
$ {\displaystyle
\;\!
\mbox{\bf v}(\cdot,s) =
\mbox{\boldmath $Q$}(\cdot,s)
} $.
We need only
be able to estimate
$ {\displaystyle
\;\!
\|\, \hat{\mbox{\boldmath $Q$}}(\cdot,s) \,
\|_{\mbox{}_{\scriptstyle \infty}}
\!\;\!
} $:
because
$ {\displaystyle
\;\!
\mathbb{F}\:\![\,\nabla \!\:\!
P(\cdot,s)\,]\:\!(k) \;\!=\;\;\!
\stackrel{\mbox{\tiny o}}{\mbox{\sf \i$\!\;\!$\i}} \!\:\!
\hat{p}(k,s) \;\!k
\,
} $
and
$ {\displaystyle
\:\!
\mbox{\small $\sum$}_{j\,=\,1}^{3} k_{j}
\hat{Q}_{j}(k,s)
=\;\!0
} $
(since
$ \nabla \!\:\!\cdot \mbox{\boldmath $Q$}(\cdot,s) =\;\! 0 $),
the vectors
$ {\displaystyle
\;\!
\mathbb{F}\:\![\,\nabla \!\:\!
P(\cdot,s)\,]\:\!(k)
} $
and
$ {\displaystyle
\:\!
\hat{\mbox{\boldmath $Q$}}(k,s)
\;\!
} $
are \linebreak
orthogonal in $\mathbb{C}^{3}\!$,
for every $k \in \mathbb{R}^{3}\!$.
Recalling
from (2.4)
that
$ {\displaystyle
\;\!
\hat{\mbox{\boldmath $Q$}}(k,s)
+\:\!
\mathbb{F}\:\![\,\nabla \!\:\!
P(\cdot,s)\;\!]\:\!(k)
=
} $ \linebreak
$ {\displaystyle
-\:
\mathbb{F}\;\![\, \mbox{\boldmath $u$}(\cdot,s) \!\;\!\cdot \!\;\!\nabla
\mbox{\boldmath $u$}(\cdot,s)\;\!]\:\!(k)
} $,
this gives \\
\mbox{} \vspace{-0.600cm} \\
\begin{equation}
\tag{2.8}
|\,\hat{\mbox{\boldmath $Q$}}(k,s) \,|_{\mbox{}_{2}}
\,\leq\;\,\!
|\: \mathbb{F}\;\![\, \mbox{\boldmath $u$}(\cdot,s)\!\;\!\cdot \!\;\!\nabla
\mbox{\boldmath $u$}\;\!(\cdot,s)\;\!]\:\!(k)
\:|_{\mbox{}_{2}}
\end{equation}
\mbox{} \vspace{-0.275cm} \\
for all $\;\! k \in \mathbb{R}^{3}\!$,
so that
we get \\
\mbox{} \vspace{-0.600cm} \\
\begin{equation}
\tag{2.9}
\|\:
\hat{\mbox{\boldmath $Q$}}(\cdot,s)
\,\|_{\mbox{}_{\scriptstyle \infty}}
\;\!\leq\;
\|\;
\mathbb{F}\;\![\, \mbox{\boldmath $u$}
\!\;\!\cdot \!\;\!\nabla
\mbox{\boldmath $u$}\,]\:\!(\cdot,s)
\:\|_{\mbox{}_{\scriptstyle \infty}}.
\end{equation}
\mbox{} \vspace{-0.200cm} \\
Now,
we have,
for each $ 1 \leq i \leq 3 $, \\
\mbox{} \vspace{-0.750cm} \\
\begin{equation}
\notag
\begin{split}
|\; \mathbb{F}\;\![\, \mbox{\boldmath $u$}(\cdot,s)
\!\;\!\cdot \!\;\!\nabla
\;\!u_{{\scriptstyle i}}(\cdot,s)\,]\:\!(k) \:|
\;\,&\leq\;
\sum_{j\,=\,1}^{3} \;\!
|\; \mathbb{F}\;\![\, u_{{\scriptstyle j}}(\cdot,s) \;\!
D_{{\scriptstyle \!\;\!j}} \:\! u_{{\scriptstyle i}}(\cdot,s)
\,]\:\!
(k) \:| \\
&\leq\;
(\:\!2\:\!\pi)^{-\,3/2} \;\!
\sum_{j\,=\,1}^{3} \;\!
\|\, u_{{\scriptstyle j}}(\cdot,s) \;\!
D_{{\scriptstyle \!\;\!j}} \:\! u_{{\scriptstyle i}}(\cdot,s) \,
\|_{\mbox{}_{\scriptstyle L^{1}(\mathbb{R}^{3})}} \\
&\leq\;
(\:\!2\:\!\pi)^{-\,3/2} \,
\|\, \mbox{\boldmath $u$}(\cdot,s) \,
\|_{\mbox{}_{\scriptstyle L^{2}(\mathbb{R}^{3})}}
\;\!
\|\, \nabla \!\;\! u_{{\scriptstyle i}}(\cdot,s) \,
\|_{\mbox{}_{\scriptstyle L^{2}(\mathbb{R}^{3})}}
\!\;\!,
\end{split}
\end{equation}
\mbox{} \vspace{+0.125cm} \\
by the Cauchy-Schwarz
inequality.
$\!$($\;\!$Here,
as before,
$ D_{\!\;\!j} \!\;\!=\,\! \partial/\partial x_{\!\;\!j} \!\;\!$.)
This gives \\
\mbox{} \vspace{-0.550cm} \\
\begin{equation}
\tag{2.10}
\|\; \mathbb{F}\;\![\, \mbox{\boldmath $u$} \!\;\!\cdot \!\;\!\nabla
\;\!\mbox{\boldmath $u$}\,]\:\!(\cdot,s) \:\|_{\mbox{}_{\scriptstyle \infty}}
\:\leq\;
(\:\!2\:\!\pi)^{-\,3/2} \,
\|\, \mbox{\boldmath $u$}(\cdot,s) \,
\|_{\mbox{}_{\scriptstyle L^{2}(\mathbb{R}^{3})}}
\;\!
\|\, D \mbox{\boldmath $u$}(\cdot,s) \,
\|_{\mbox{}_{\scriptstyle L^{2}(\mathbb{R}^{3})}}
\!\:\!.
\end{equation}
\mbox{} \vspace{-0.200cm} \\
From (2.7), (2.9) and (2.10),
one gets (2.5$a$),
which shows the first part
of Theorem 2.1. \linebreak
\mbox{} \vspace{-0.750cm} \\

The proof of (2.5$b$) follows
in a similar way,
using (2.8)
and the elementary estimate \linebreak
\mbox{} \vspace{-0.530cm} \\
$ {\displaystyle
\|\: e^{\:\!\Delta \mbox{\footnotesize $\tau$}}
\!\;\!\mbox{u} \:
\|_{\mbox{}_{\scriptstyle L^{\infty}(\mathbb{R}^{3})}}
\!\:\!\leq\;\!
K \;\! \tau^{-\,3/4} \,
\|\: \mbox{u} \:
\|_{\mbox{}_{\scriptstyle L^{2}(\mathbb{R}^{3})}}
\!\;\!
} $
for the heat semigroup,
where
$ \tau > 0 $ is arbitrary,
and \linebreak
\mbox{} \vspace{-0.530cm} \\
$ K \!\:\!=\:\! (\:\! 8 \:\!\pi )^{-\,3/4} \!\;\!$.
$\!$This gives,
for any
$ \:\!s > 0 \;\! $
with
both
$ {\displaystyle
\;\!
\|\, \mbox{\boldmath $u$}(\cdot,s) \,
\|_{\mbox{}_{\scriptstyle \infty}}
\!\:\!
} $
and
$ {\displaystyle
\:\!
\|\, D \mbox{\boldmath $u$}(\cdot,s) \,
\|_{\mbox{}_{\scriptstyle L^{2}(\mathbb{R}^{3})}}
\!
} $
finite, \\
\mbox{} \vspace{-0.050cm} \\
\mbox{} \hspace{+1.250cm}
$ {\displaystyle
\|\: e^{\:\!\mbox{\scriptsize $\Delta$} \mbox{\scriptsize $(t - s)$}}
\;\! \mbox{\boldmath $Q$}(\cdot,s) \,
\|_{\mbox{}_{\scriptstyle \infty}}
\,\leq\;
K \;\!
(\:\!t - s )^{-\,3/4} \,
\|\, \mbox{\boldmath $Q$}(\cdot,s) \,
\|_{\mbox{}_{\scriptstyle L^{2}(\mathbb{R}^{3})}}
} $ \\
\mbox{} \vspace{-0.150cm} \\
\mbox{} \hspace{+4.725cm}
$ {\displaystyle
\leq\;
K \;\!
(\:\!t - s )^{-\,3/4} \,
\|\, \mbox{\boldmath $u$}(\cdot,s)
\!\:\!\cdot\!\:\! \nabla \:\! \mbox{\boldmath $u$}(\cdot,s) \,
\|_{\mbox{}_{\scriptstyle L^{2}(\mathbb{R}^{3})}}
} $
\mbox{} \hfill [$\,$by (2.8)$\,$] \\
\mbox{} \vspace{-0.150cm} \\
\mbox{} \hspace{+4.725cm}
$ {\displaystyle
\leq\;
K \;\!
(\:\!t - s )^{-\,3/4} \,
\|\, \mbox{\boldmath $u$}(\cdot,s) \,
\|_{\mbox{}_{\scriptstyle \infty}} \:\!
\|\, D\mbox{\boldmath $u$}(\cdot,s) \,
\|_{\mbox{}_{\scriptstyle L^{2}(\mathbb{R}^{3})}}
} $ \\
%
%
\mbox{} \vspace{-0.050cm} \\
for all $\;\!t > s $,
using
Parseval's identity (twice)
and the norm definitions
(1.16), (1.17).
}
\mbox{} \hfill $\Box$ \\
%
%
\mbox{} \vspace{-0.750cm} \\

Let us notice that,
applying the argument above
to solutions
of the regularized Navier-Stokes equations (2.1),
we obtain,
in a completely similar way, \\
\mbox{} \vspace{-0.570cm} \\
\begin{equation}
\tag{2.11$a$}
\|\: e^{\:\!\mbox{\scriptsize $\Delta$} \mbox{\scriptsize $(t - s)$}} \;\!
\mbox{\boldmath $Q$}_{\mbox{}_{\scriptstyle \!\delta}}\!\;\!(\cdot,s) \,
\|_{\mbox{}_{\scriptstyle L^{2}(\mathbb{R}^{3})}}
\;\!\leq\,
K \;\!
(t - s)^{-\,3/4} \,
\|\, \mbox{\boldmath $u$}_{\mbox{}_{\scriptstyle \!\:\!\delta}}
\!\;\!(\cdot,s) \,
\|_{\mbox{}_{\scriptstyle L^{2}(\mathbb{R}^{3})}}
\;\!
\|\, D \mbox{\boldmath $u$}_{\mbox{}_{\scriptstyle \!\:\!\delta}}
\!\;\!(\cdot,s) \,
\|_{\mbox{}_{\scriptstyle L^{2}(\mathbb{R}^{3})}}
\end{equation}
\mbox{} \vspace{-0.300cm} \\
and \\
\mbox{} \vspace{-0.800cm} \\
\begin{equation}
\tag{2.11$b$}
\|\: e^{\:\!\mbox{\scriptsize $\Delta$} \mbox{\scriptsize $(t - s)$}} \;\!
\mbox{\boldmath $Q$}_{\mbox{}_{\scriptstyle \!\delta}}\!\;\!(\cdot,s) \,
\|_{\mbox{}_{\scriptstyle \infty}}
\;\!\leq\,
K \;\!
(t - s)^{-\,3/4} \,
\|\, \mbox{\boldmath $u$}_{\mbox{}_{\scriptstyle \!\:\!\delta}}
\!\;\!(\cdot,s) \,
\|_{\mbox{}_{\scriptstyle \infty}}
\;\!
\|\, D \mbox{\boldmath $u$}_{\mbox{}_{\scriptstyle \!\:\!\delta}}
\!\;\!(\cdot,s) \,
\|_{\mbox{}_{\scriptstyle L^{2}(\mathbb{R}^{3})}}
\end{equation}
\mbox{} \vspace{-0.250cm} \\
for all $\;\!t > s > 0 $,
where
$ K \!\;\!=\:\! (\:\!8 \:\!\pi)^{-\,3/4} \!$,
as before,
and \\
\mbox{} \vspace{-0.600cm} \\
\begin{equation}
\tag{2.12}
\mbox{\boldmath $Q$}_{\mbox{}_{\scriptstyle \!\delta}}
\!\:\!(\cdot,s)
\;=\; -\:
\bar{\mbox{\boldmath $u$}}_{\mbox{}_{\scriptstyle \!\:\!\delta}}\!\:\!(\cdot,s)
\!\;\!\cdot \!\;\!\nabla
\mbox{\boldmath $u$}_{\mbox{}_{\scriptstyle \!\:\!\delta}}\!\:\!(\cdot,s)
\,-\;\!
\nabla p_{\mbox{}_{\scriptstyle \!\delta}}\!\;\!(\cdot,s).
\end{equation}
\mbox{} \vspace{-0.250cm} \\
The estimate (2.11$a$)
is particularly useful,
since
%
%
the regularized solutions
$ {\displaystyle
\:\!
\mbox{\boldmath $u$}_{\mbox{}_{\scriptstyle \!\:\!\delta}}\!\:\!(\cdot,t)
\:\!
} $
given in~(2.1)
satisfy
the energy inequality \\
\mbox{} \vspace{-0.650cm} \\
\begin{equation}
\tag{2.13}
\|\,\mbox{\boldmath $u$}_{\mbox{}_{\scriptstyle \!\:\!\delta}}
\!\:\!(\cdot,t) \,
\|_{\mbox{}_{\scriptstyle L^{2}(\mathbb{R}^{3})}}^{\:\!2}
+\;
2 \!
\int_{0}^{\;\!\mbox{\mbox{\footnotesize $t$}}} \!\!\;\!
\|\, D \mbox{\boldmath $u$}_{\mbox{}_{\scriptstyle \!\:\!\delta}}
\!\:\!(\cdot,s) \,
\|_{\mbox{}_{\scriptstyle L^{2}(\mathbb{R}^{3})}}^{\:\!2}
\;\!ds
\;\;\!\leq\:
\|\, \mbox{\boldmath $u$}_{0} \,
\|_{\mbox{}_{\scriptstyle L^{2}(\mathbb{R}^{3})}}^{\:\!2}
\end{equation}
\mbox{} \vspace{-0.100cm} \\
for all
$\;\! t > 0 $
(and
$ \delta > 0 $
arbitrary),
from which
$ {\displaystyle
\;\!
\|\,\mbox{\boldmath $u$}_{\mbox{}_{\scriptstyle \!\:\!\delta}}
\!\:\!(\cdot,t) \,
\|_{\mbox{}_{\scriptstyle L^{2}(\mathbb{R}^{3})}}
\!\;\!
} $,
$
\int_{0}^{\;\!\mbox{\mbox{\footnotesize $t$}}} \!
\|\, D \mbox{\boldmath $u$}_{\mbox{}_{\scriptstyle \!\:\!\delta}}
\!\:\!(\cdot,s) \,
\|_{\mbox{}_{\scriptstyle L^{2}(\mathbb{R}^{3})}}^{\:\!2}
\!\;\!ds
$ \linebreak
\mbox{} \vspace{-0.550cm} \\
can be bounded
independently of $ \:\!\delta > 0 $.
This will be used in Theorems 2.2 and 2.3 below
to show that
the particular
value of
$\;\!t_0 \!\;\!\geq 0 \;\!$
chosen in
defining the heat flow approximations
%
%
(1.7)
is not relevant
in regard to
the properties
(1.8)$\;\!-\;\!$(1.10). \\
\nl
%
%
%
%
{\bf Theorem 2.2.}
\textit{%
Let
$ \;\!\mbox{\boldmath $u$}(\cdot,t) $,
$ t > 0 $,
be any particular
Leray-Hopf's solution to $\;\!(1.1)$.
Given
any pair of
initial values
$\;\! \tilde{t}_0 \!\;\!> t_0 \!\;\!\geq 0 $,
one has
} \\
\mbox{} \vspace{-0.575cm} \\
\begin{equation}
\tag{2.14}
\|\, \mbox{\boldmath $v$}(\cdot,t) \,-\,
\tilde{\mbox{\boldmath $v$}}(\cdot,t) \,
\|_{\mbox{}_{\scriptstyle L^{2}(\mathbb{R}^{3})}}
\,\leq\;
\frac{\,\mbox{\small $K$}}{\mbox{\small $\sqrt{\;\!2\;}\;$}}
\: \|\, \mbox{\boldmath $u$}_0 \,
\|_{\mbox{}_{\scriptstyle L^{2}(\mathbb{R}^{3})}}^{\:\!2}
(\:\! \tilde{t}_{0} \!\;\!-\:\! t_0 )^{1/2}
\,
(\:\! t -\:\! \tilde{t}_0 )^{-\,3/4}
\end{equation}
\mbox{} \vspace{-0.075cm} \\
\textit{%
for all $\,t > \tilde{t}_0 $,
where
$ {\displaystyle
\,
\mbox{\boldmath $v$}(\cdot,t)
\:\!=\;\!
e^{\:\!\mbox{\scriptsize $\Delta$} (\:\!\mbox{\footnotesize $t$}
\;\!-\,\mbox{\footnotesize $t_0$})}
\:\!\mbox{\boldmath $u$}(\cdot,t_0)
} $,
$ {\displaystyle
\:\!
\tilde{\mbox{\boldmath $v$}}(\cdot,t)
\:\!=\;\!
e^{\:\!\mbox{\scriptsize $\Delta$} (\:\!\mbox{\footnotesize $t$}
\;\!-\,\mbox{\footnotesize $\tilde{t}_0$})}
\:\!\mbox{\boldmath $u$}(\cdot,\tilde{t}_0)
} $
are the corresponding
heat flows
associated with
$ \;\! t_0 $, $ \tilde{t}_0 $,
respectively,
and
$ {\displaystyle
\;\!
\mbox{\small $K$}
\!\:\!=\:\!
(\:\!8 \:\!\pi)^{-\,3/4}
\!\:\!
} $.
} \\
%
%
\nl
\mbox{} \vspace{-0.450cm} \\
{\small
{\bf Proof:}
We start by writing
$ \:\!\mbox{\boldmath $v$}(\cdot,t) \:\!$
as \\
\mbox{} \vspace{-0.500cm} \\
\begin{equation}
\notag
\mbox{\boldmath $v$}(\cdot,t)
\;=\;\:\!
e^{\;\!\Delta \:\!(\;\!\mbox{\footnotesize $t$} \;\!-\,
\mbox{\footnotesize $t_0$})} \;\!
[\,\mbox{\boldmath $u$}(\cdot,t_0) -
\mbox{\boldmath $u$}_{\mbox{}_{\scriptstyle \!\:\!\delta}}
\!\;\!(\cdot,t_0) \,]
\:+\:
e^{\;\!\Delta \:\!(\;\!\mbox{\footnotesize $t$} \;\!-\,
\mbox{\footnotesize $t_0$})} \;\!
\mbox{\boldmath $u$}_{\mbox{}_{\scriptstyle \!\:\!\delta}}
\!\;\!(\cdot,t_0),
\qquad
t > t_0,
\end{equation}
\mbox{} \vspace{-0.200cm} \\
with
$ {\displaystyle
\;\!
\mbox{\boldmath $u$}_{\mbox{}_{\scriptstyle \!\:\!\delta}}
\!\;\!(\cdot,t)
\;\!
} $
given in (2.1),
$ \delta > 0 $.
Because \\
\mbox{} \vspace{-0.600cm} \\
\begin{equation}
\notag
\mbox{\boldmath $u$}_{\mbox{}_{\scriptstyle \!\:\!\delta}}
\!\;\!(\cdot,t_0)
\;=\;\:\!
e^{\;\!\Delta \:\!\mbox{\footnotesize $t_0$}}
\;\!\bar{\mbox{\boldmath $u$}}_{0, \,\delta}
\:+\!\:\!
\int_{\:\!\mbox{\footnotesize $0$}}^{\;\!\mbox{\footnotesize $t_0$}}
\!\!
e^{\;\!\Delta \:\!(\;\!\mbox{\footnotesize $t_0$} \;\!-\,
\mbox{\footnotesize $s$})} \;\!
\mbox{\boldmath $Q$}_{\mbox{}_{\scriptstyle \!\delta}}
\!\;\!(\cdot,s)
\,ds,
\end{equation}
\mbox{} \vspace{-0.125cm} \\
where
$ {\displaystyle
\;\!
\bar{\mbox{\boldmath $u$}}_{\mbox{}_{\scriptstyle \!\;\!0, \,\delta}}
\!\;\!=\;\!
G_{\mbox{}_{\scriptstyle \!\delta}}
\!\:\!\ast
\mbox{\boldmath $u$}_{0}
} $,
$ {\displaystyle
\:\!
\mbox{\boldmath $Q$}_{\mbox{}_{\scriptstyle \!\delta}}
\!\:\!(\cdot,s)
= -\,
\bar{\mbox{\boldmath $u$}}_{\mbox{}_{\scriptstyle \!\:\!\delta}}\!\:\!(\cdot,s)
\!\;\!\cdot \!\;\!\nabla
\mbox{\boldmath $u$}_{\mbox{}_{\scriptstyle \!\:\!\delta}}\!\:\!(\cdot,s)
\;\!-\:\!
\nabla p_{\mbox{}_{\scriptstyle \!\delta}}\!\;\!(\cdot,s)
} $,
$\;\!$cf.$\;$(2.1$b$) and (2.12) above,
we get \\
\mbox{} \vspace{-0.800cm} \\
\begin{equation}
\notag
\mbox{\boldmath $v$}(\cdot,t)
\;=\;\:\!
e^{\;\!\Delta \:\!(\;\!\mbox{\footnotesize $t$} \;\!-\,
\mbox{\footnotesize $t_0$})} \;\!
[\,\mbox{\boldmath $u$}(\cdot,t_0) -
\:\!\mbox{\boldmath $u$}_{\mbox{}_{\scriptstyle \!\:\!\delta}}
\!\;\!(\cdot,t_0) \,]
\;+\;
e^{\;\!\Delta \mbox{\footnotesize $t$}} \:\!
\bar{\mbox{\boldmath $u$}}_{0,\,\delta}
\;+
\int_{\:\!\mbox{\footnotesize $0$}}^{\;\!\mbox{\footnotesize $t_0$}}
\!\!
e^{\;\!\Delta \:\!(\;\!\mbox{\footnotesize $t$} \;\!-\,
\mbox{\footnotesize $s$})} \;\!
\mbox{\boldmath $Q$}_{\mbox{}_{\scriptstyle \!\delta}}
\!\:\!(\cdot,s)
\; ds,
\end{equation}
\mbox{} \vspace{-0.225cm} \\
for $ \;\!t > t_0 $.
Similarly,
we have,
for $\;\! t > \tilde{t}_0 $: \\
\mbox{} \vspace{-0.750cm} \\
\begin{equation}
\notag
\tilde{\mbox{\boldmath $v$}}(\cdot,t)
\;=\;\:\!
e^{\;\!\Delta \:\!(\;\!\mbox{\footnotesize $t$} \;\!-\,
\mbox{\footnotesize $\tilde{t}_0$})} \;\!
[\,\mbox{\boldmath $u$}(\cdot,\tilde{t}_0) -
\:\!\mbox{\boldmath $u$}_{\mbox{}_{\scriptstyle \!\:\!\delta}}
\!\;\!(\cdot,\tilde{t}_0) \,]
\;+\;
e^{\;\!\Delta \mbox{\footnotesize $t$}} \:\!
\bar{\mbox{\boldmath $u$}}_{0,\,\delta}
\;+
\int_{\:\!\mbox{\footnotesize $0$}}^{\;\!\mbox{\footnotesize $\tilde{t}_0$}}
\!\!
e^{\;\!\Delta \:\!(\;\!\mbox{\footnotesize $t$} \;\!-\,
\mbox{\footnotesize $s$})} \;\!
\mbox{\boldmath $Q$}_{\mbox{}_{\scriptstyle \!\delta}}
\!\:\!(\cdot,s)
\; ds.
\end{equation}
\mbox{} \vspace{-0.150cm} \\
Hence,
we obtain,
for the difference
$ {\displaystyle
\;\!
\mbox{\boldmath $v$}(\cdot,t)
\;\!-\;\!
\tilde{\mbox{\boldmath $v$}}(\cdot,t)
} $,
at any $\:\! t > \tilde{t}_0 $,
the identity \\
\mbox{} \vspace{-0.200cm} \\
\mbox{} \hspace{+0.025cm}
$ {\displaystyle
\tilde{\mbox{\boldmath $v$}}(\cdot,t)
\,-\,
\mbox{\boldmath $v$}(\cdot,t)
\;\;\!=\;\;\;\!
e^{\;\!\Delta \:\!(\;\!\mbox{\footnotesize $t$}
\;\!-\,
\mbox{\footnotesize $\tilde{t}_0$})} \;\!
[\,\mbox{\boldmath $u$}(\cdot,\tilde{t}_0) -
\:\!\mbox{\boldmath $u$}_{\mbox{}_{\scriptstyle \!\:\!\delta}}
\!\;\!(\cdot,\tilde{t}_0) \,]
\;-\;
e^{\;\!\Delta \:\!(\;\!\mbox{\footnotesize $t$} \;\!-\,
\mbox{\footnotesize $t_0$})} \;\!
[\,\mbox{\boldmath $u$}(\cdot,t_0) -
\:\!\mbox{\boldmath $u$}_{\mbox{}_{\scriptstyle \!\:\!\delta}}
\!\;\!(\cdot,t_0) \,]
} $ \\
\mbox{} \vspace{-0.250cm} \\
\mbox{} \hspace{+8.400cm}
$ {\displaystyle
+\:\!
\int_{\:\!\mbox{\footnotesize $t_0$}}
    ^{\;\!\mbox{\footnotesize $\tilde{t}_0$}}
\!\!
e^{\;\!\Delta \:\!(\;\!\mbox{\footnotesize $t$} \;\!-\,
\mbox{\footnotesize $s$})} \;\!
\mbox{\boldmath $Q$}_{\mbox{}_{\scriptstyle \!\delta}}
\!\:\!(\cdot,s)
\; ds
} $.
\mbox{} \hfill (2.15) \\
\mbox{} \vspace{+0.200cm} \\
Therefore,
given any
$ {\displaystyle
\;\!\mathbb{K} \subset \mathbb{R}^{3}
\!\:\!
} $
compact,
we get,
for each
$ \;\! t > \tilde{t}_0 $,
$ \;\! \delta > 0 $: \\
\mbox{} \vspace{-0.750cm} \\
\begin{equation}
\notag
\begin{split}
\|\, \tilde{\mbox{\boldmath $v$}}(\cdot,t)
\,-\, \mbox{\boldmath $v$}(\cdot,t) \,
\|_{\mbox{}_{\scriptstyle L^{2}(\mathbb{K})}}
\,&\leq\;
J_{\mbox{}_{\scriptstyle \!\:\!\delta}}\!\:\!(t)
\;+\:\!
\int_{\:\!\mbox{\footnotesize $t_0$}}
    ^{\;\!\mbox{\footnotesize $\tilde{t}_0$}}
\!\!
\|\:
e^{\;\!\Delta \:\!(\;\!\mbox{\footnotesize $t$} \;\!-\,
\mbox{\footnotesize $s$})} \;\!
\mbox{\boldmath $Q$}_{\mbox{}_{\scriptstyle \!\delta}}
\!\:\!(\cdot,s)
\:\|_{\mbox{}_{\scriptstyle L^{2}(\mathbb{K})}}
\, ds \\
&\leq\;
J_{\mbox{}_{\scriptstyle \!\:\!\delta}}\!\;\!(t)
\;+\,
K \!\!\:\!
\int_{\:\!\mbox{\footnotesize $t_0$}}
    ^{\;\!\mbox{\footnotesize $\tilde{t}_0$}}
\!\!
(\:\! t - s \:\!)^{-\,3/4} \:
\|\,\mbox{\boldmath $u$}_{\mbox{}_{\scriptstyle \!\:\!\delta}}
\!\:\!(\cdot,s) \,
\|_{\mbox{}_{\scriptstyle L^{2}(\mathbb{R}^{3})}}
\:\!
\|\, D \mbox{\boldmath $u$}_{\mbox{}_{\scriptstyle \!\:\!\delta}}
\!\:\!(\cdot,s) \,
\|_{\mbox{}_{\scriptstyle L^{2}(\mathbb{R}^{3})}}
\;\! ds \\
&\leq\;
J_{\mbox{}_{\scriptstyle \!\:\!\delta}}\!\:\!(t)
\;+\;
\frac{\;K}{\sqrt{\,2\;}\;} \,
(\;\! \tilde{t}_0 \!\;\!-\;\! t_0 )^{\mbox{}^{\scriptstyle
\!\:\! \frac{\scriptstyle 1}{\scriptstyle 2} }}
\:\!
\|\, \mbox{\boldmath $u$}_0 \,
\|_{\mbox{}_{\scriptstyle L^{2}(\mathbb{R}^{3})}}
  ^{\mbox{}^{\scriptstyle \:\! 2}}
\;\!
(\;\! t - \tilde{t}_0 )^{\mbox{}^{\scriptstyle
\!\!\!\:\! -\, \frac{\scriptstyle 3}{\scriptstyle 4} }}
\end{split}
\end{equation}
\mbox{} \vspace{+0.075cm} \\
by (2.11$a$), (2.13),
where
$ \:\!K \!\:\!=\:\!(\:\!8\:\!\pi )^{-\,3/4} \!\;\!$
and \\
\mbox{} \vspace{+0.000cm} \\
\mbox{} \hspace{+0.100cm}
$ {\displaystyle
J_{\mbox{}_{\scriptstyle \!\:\!\delta}}\!\:\!(t)
\;=\;
\|\: e^{\;\!\Delta \:\!(\;\!\mbox{\footnotesize $t$}
\;\!-\, \mbox{\footnotesize $\tilde{t}_0$})} \;\!
[\,\mbox{\boldmath $u$}(\cdot,\tilde{t}_0) -
\:\!\mbox{\boldmath $u$}_{\mbox{}_{\scriptstyle \!\:\!\delta}}
\!\;\!(\cdot,\tilde{t}_0) \,]
\:\|_{\mbox{}_{\scriptstyle L^{2}(\mathbb{K})}}
\:\!+\;
\|\: e^{\;\!\Delta \:\!(\;\!\mbox{\footnotesize $t$}
\;\!-\, \mbox{\footnotesize $t_0$})} \;\!
[\,\mbox{\boldmath $u$}(\cdot,t_0) -
\:\!\mbox{\boldmath $u$}_{\mbox{}_{\scriptstyle \!\:\!\delta}}
\!\;\!(\cdot,t_0) \,]
\: \|_{\mbox{}_{\scriptstyle L^{2}(\mathbb{K})}}
\!\:\!
} $. \\
\mbox{} \vspace{+0.050cm} \\
Taking
$ \;\!\delta = \delta^{\prime} \!\rightarrow  0 \;\! $
according to (2.2),
we get
$ J_{\mbox{}_{\scriptstyle \!\:\!\delta}}\!\;\!(t) \rightarrow 0 $,
since,
by Lebesgue's Dominated Convergence Theorem
and (2.2),
we have,
for any
$ \sigma \!\:\!, \, \tau > 0 \:\!$: \\
\mbox{} \vspace{-0.550cm} \\
\begin{equation}
\notag
\|\: e^{\;\!\Delta \mbox{\footnotesize $\tau$}}
[\,\mbox{\boldmath $u$}(\cdot,\sigma) -
\:\!\mbox{\boldmath $u$}_{\mbox{}_{\scriptstyle \!\:\!\delta^{\prime}}}
\!\;\!(\cdot,\sigma) \,]
\:\|_{\mbox{}_{\scriptstyle L^{2}(\mathbb{K})}}
\;\!\rightarrow\; 0
\qquad
\mbox{as }
\;\; \delta^{\prime} \!\rightarrow 0,
\end{equation}
\mbox{} \vspace{-0.190cm} \\
recalling that
$ \mathbb{K} $
has finite measure.
Hence,
we obtain \\
\mbox{} \vspace{-0.550cm} \\
\begin{equation}
\notag
\|\, \tilde{\mbox{\boldmath $v$}}(\cdot,t)
\,-\, \mbox{\boldmath $v$}(\cdot,t) \,
\|_{\mbox{}_{\scriptstyle L^{2}(\mathbb{K})}}
\leq\:
\frac{\;K}{\sqrt{\;\!2\;}\;} \,
(\;\! \tilde{t}_0 \!\;\!-\;\! t_0 )^{\!\;\!1/2}
\,
\|\, \mbox{\boldmath $u$}_0 \,
\|_{\mbox{}_{\scriptstyle L^{2}(\mathbb{R}^{3})}}
  ^{\scriptstyle \:\!2}
\;\!
(\;\! t - \tilde{t}_0 )^{-\,3/4}
\end{equation}
\mbox{} \vspace{-0.170cm} \\
for each
$ \;\! t > \tilde{t}_0 $,
and for
{\em any\/}
compact set
$ \;\!\mathbb{K} \subset \mathbb{R}^{3} \!\;\!$.
This is clearly equivalent to (2.14).
}
\mbox{} \hfill $\Box$ \\
%
%
\mbox{} \vspace{-0.675cm} \\

Theorem 2.2 greatly simplifies
the derivation of
the asymptotic property (1.8).
For similar reasons,
our proof of (1.9)
requires the supnorm version
of (2.14) above,
which is given
in the next result. \\
\mbox{} \vspace{-0.050cm} \\
%
%
%
%
{\bf Theorem 2.3.}
\textit{%
Let
$ \;\!\mbox{\boldmath $u$}(\cdot,t) $,
$ t > 0 $,
be any particular
Leray-Hopf's solution to $\;\!(1.1)$.
Given
any pair of
initial values
$\;\! \tilde{t}_0 \!\;\!> t_0 \!\;\!\geq 0 $,
one has
} \\
\mbox{} \vspace{-0.620cm} \\
\begin{equation}
\tag{2.16}
\|\, \mbox{\boldmath $v$}(\cdot,t) \,-\,
\tilde{\mbox{\boldmath $v$}}(\cdot,t) \,
\|_{\mbox{}_{\scriptstyle L^{\infty}(\mathbb{R}^{3})}}
\;\!\leq\;
\frac{\,\mbox{\small $\Gamma$}}{\mbox{\small $\sqrt{\;\!2\;}\;$}}
\: \|\, \mbox{\boldmath $u$}_0 \,
\|_{\mbox{}_{\scriptstyle L^{2}(\mathbb{R}^{3})}}^{\:\!2}
(\:\! \tilde{t}_{0} \!\;\!-\:\! t_0 )^{1/2}
\;\!
(\:\! t -\:\! \tilde{t}_0 )^{-\,3/2}
\end{equation}
\mbox{} \vspace{-0.030cm} \\
\textit{%
for all $\,t > \tilde{t}_0 $,
where
$ {\displaystyle
\,
\mbox{\boldmath $v$}(\cdot,t)
\:\!=\;\!
e^{\:\!\mbox{\scriptsize $\Delta$} (\:\!\mbox{\footnotesize $t$}
\;\!-\,\mbox{\footnotesize $t_0$})}
\:\!\mbox{\boldmath $u$}(\cdot,t_0)
} $,
$ {\displaystyle
\:\!
\tilde{\mbox{\boldmath $v$}}(\cdot,t)
\:\!=\;\!
e^{\:\!\mbox{\scriptsize $\Delta$} (\:\!\mbox{\footnotesize $t$}
\;\!-\,\mbox{\footnotesize $\tilde{t}_0$})}
\:\!\mbox{\boldmath $u$}(\cdot,\tilde{t}_0)
} $
are the corresponding
heat flows
associated with
$ \;\! t_0 $, $ \tilde{t}_0 $,
respectively,
and
$ {\displaystyle
\,
\mbox{\small $\Gamma$}
=\:\!
(\:\!4 \:\!\pi)^{-\,3/2}
\!\:\!
} $.
}

%
%
\nl
%
{\small
{\bf Proof:}
Taking
$ {\displaystyle
\;\!\mathbb{K} \subset \mathbb{R}^{3}
\!\:\!
} $
compact
and
$ \:\! 2 < q < \infty \:\!$
arbitrary,
we get,
for each
$ \;\! t > \tilde{t}_0 $,
$ \;\! \delta > 0 $,
recalling (2.12), (2.15): \\
\mbox{} \vspace{-0.850cm} \\
\begin{equation}
\notag
\begin{split}
\|\, \tilde{\mbox{\boldmath $v$}}(\cdot,t)
\,-\, \mbox{\boldmath $v$}(\cdot,t) \,
\|_{\mbox{}_{\scriptstyle L^{q}(\mathbb{K})}}
\,&\leq\;
J_{\mbox{}_{\scriptstyle \!\:\!\delta, \;\!q}}\!\:\!(t)
\;+\:\!
\int_{\:\!\mbox{\footnotesize $t_0$}}
    ^{\;\!\mbox{\footnotesize $\tilde{t}_0$}}
\!\!
\|\:
e^{\;\!\Delta \:\!(\;\!\mbox{\footnotesize $t$} \;\!-\,
\mbox{\footnotesize $s$})} \;\!
\mbox{\boldmath $Q$}_{\mbox{}_{\scriptstyle \!\delta}}
\!\:\!(\cdot,s)
\:\|_{\mbox{}_{\scriptstyle L^{q}(\mathbb{R}^{3})}}
\, ds \\
&\leq\;
J_{\mbox{}_{\scriptstyle \!\:\!\delta, \;\!q}}\!\:\!(t)
\;+
\int_{\:\!\mbox{\footnotesize $t_0$}}
    ^{\;\!\mbox{\footnotesize $\tilde{t}_0$}}
\!
\bigl[\;\! 4 \:\!\pi \;\!(\:\!t - s ) \:\!\bigr]^{\mbox{}^{\scriptstyle \!\!
-\, \frac{\scriptstyle 3}{\scriptstyle 4}
\bigl( 1 \,-\, \frac{\scriptstyle 2}{\scriptstyle q} \bigr) }}
\!
\|\:
e^{\mbox{}^{\scriptstyle \!\:\!\frac{1}{2} \:\! \Delta \:\!(\;\!\mbox{\footnotesize $t$} \;\!-\,
\mbox{\footnotesize $s$})} } \!\:\!
\mbox{\boldmath $Q$}_{\mbox{}_{\scriptstyle \!\delta}}
\!\:\!(\cdot,s)
\:\|_{\mbox{}_{\scriptstyle L^{2}(\mathbb{R}^{3})}}
\, ds \\
&\leq\;
J_{\mbox{}_{\scriptstyle \!\:\!\delta, \;\!q}}\!\;\!(t)
\:+\:
\gamma_{q} \!\!\:\!
\int_{\:\!\mbox{\footnotesize $t_0$}}
    ^{\;\!\mbox{\footnotesize $\tilde{t}_0$}}
\!\!
(\:\! t - s \:\!)^{\mbox{}^{\scriptstyle \!\!\!
-\, \frac{\scriptstyle 3}{\scriptstyle 2}
\bigl( 1 \,-\, \frac{\scriptstyle 1}{\scriptstyle q} \bigr) }}
\!\:\!
\|\,\mbox{\boldmath $u$}_{\mbox{}_{\scriptstyle \!\:\!\delta}}
\!\:\!(\cdot,s) \,
\|_{\mbox{}_{\scriptstyle L^{2}(\mathbb{R}^{3})}}
\|\, D \mbox{\boldmath $u$}_{\mbox{}_{\scriptstyle \!\:\!\delta}}
\!\:\!(\cdot,s) \,
\|_{\mbox{}_{\scriptstyle L^{2}(\mathbb{R}^{3})}} \\
&\leq\;
J_{\mbox{}_{\scriptstyle \!\:\!\delta,\;\!q}}\!\:\!(t)
\;+\;
\frac{\:\gamma_{q}}{\sqrt{\;\!2\;}\;} \,
(\;\! \tilde{t}_0 \!\;\!-\;\! t_0 )^{\mbox{}^{\scriptstyle
\!\:\! \frac{\scriptstyle 1}{\scriptstyle 2} }}
\:\!
\|\, \mbox{\boldmath $u$}_0 \,
\|_{\mbox{}_{\scriptstyle L^{2}(\mathbb{R}^{3})}}
  ^{\mbox{}^{\scriptstyle \:\! 2}}
\;\!
(\;\! t - \tilde{t}_0 )^{\mbox{}^{\scriptstyle \!\!\!
-\, \frac{\scriptstyle 3}{\scriptstyle 2}
\bigl( 1 \,-\, \frac{\scriptstyle 1}{\scriptstyle q} \bigr) }}
\!\:\!
\end{split}
\end{equation}
\mbox{} \vspace{+0.050cm} \\
by (2.11$a$), (2.13),
where
$ {\displaystyle
\;\!
\gamma_{q} \!\;\!=
(\:\! 4 \:\! \pi)^{\mbox{}^{\scriptstyle \!\!\!
-\, \frac{\scriptstyle 3}{\scriptstyle 2}
\bigl( 1 \,-\, \frac{\scriptstyle 1}{\scriptstyle q} \bigr) }}
\!\!
} $
and \\
\mbox{} \vspace{-0.025cm} \\
\mbox{} \hspace{+0.100cm}
$ {\displaystyle
J_{\mbox{}_{\scriptstyle \!\:\!\delta, \;\!q}}\!\:\!(t)
\:=\;
\|\: e^{\;\!\Delta \:\!(\;\!\mbox{\footnotesize $t$}
\;\!-\, \mbox{\footnotesize $\tilde{t}_0$})} \;\!
[\,\mbox{\boldmath $u$}(\cdot,\tilde{t}_0) -
\:\!\mbox{\boldmath $u$}_{\mbox{}_{\scriptstyle \!\:\!\delta}}
\!\;\!(\cdot,\tilde{t}_0) \,]
\:\|_{\mbox{}_{\scriptstyle L^{q}(\mathbb{K})}}
\:\!+\;
\|\: e^{\;\!\Delta \:\!(\;\!\mbox{\footnotesize $t$}
\;\!-\, \mbox{\footnotesize $t_0$})} \;\!
[\,\mbox{\boldmath $u$}(\cdot,t_0) -
\:\!\mbox{\boldmath $u$}_{\mbox{}_{\scriptstyle \!\:\!\delta}}
\!\;\!(\cdot,t_0) \,]
\: \|_{\mbox{}_{\scriptstyle L^{q}(\mathbb{K})}}
\!\:\!
} $. \\
\mbox{} \vspace{+0.000cm} \\
Taking
$ \;\!\delta = \delta^{\prime} \!\rightarrow  0 \;\! $
according to (2.2),
we get
$ J_{\mbox{}_{\scriptstyle \!\:\!\delta,\;\!q}}\!\;\!(t) \rightarrow 0 $,
since,
by Lebesgue's Dominated Convergence Theorem
and (2.2),
we have
$ {\displaystyle
\;\!
\|\: e^{\;\!\Delta \mbox{\footnotesize $\tau$}}
[\,\mbox{\boldmath $u$}(\cdot,\sigma) -
\:\!\mbox{\boldmath $u$}_{\mbox{}_{\scriptstyle \!\:\!\delta^{\prime}}}
\!\;\!(\cdot,\sigma) \,]
\:\|_{\mbox{}_{\scriptstyle L^{q}(\mathbb{K})}}
\;\!\rightarrow\; 0
\;\!
} $
as
$ {\displaystyle
\:\! \delta^{\prime} \!\rightarrow 0
} $,
for each
$ \:\!\sigma \!\:\!, \, \tau > 0 $.
Hence,
letting
$ \;\!\delta = \delta^{\prime} \!\rightarrow  0 $,
we obtain \\
\mbox{} \vspace{-0.750cm} \\
\begin{equation}
\notag
\|\, \tilde{\mbox{\boldmath $v$}}(\cdot,t)
\,-\, \mbox{\boldmath $v$}(\cdot,t) \,
\|_{\mbox{}_{\scriptstyle L^{q}(\mathbb{K})}}
\;\!\leq\:
\frac{\:\gamma_{q}}{\sqrt{\;\!2\;}\;} \,
(\;\! \tilde{t}_0 \!\;\!-\;\! t_0 )^{\mbox{}^{\scriptstyle
\!\:\! \frac{\scriptstyle 1}{\scriptstyle 2} }}
\:\!
\|\, \mbox{\boldmath $u$}_0 \,
\|_{\mbox{}_{\scriptstyle L^{2}(\mathbb{R}^{3})}}
  ^{\mbox{}^{\scriptstyle \:\! 2}}
(\;\! t - \tilde{t}_0 )^{\mbox{}^{\scriptstyle \!\!\!
-\, \frac{\scriptstyle 3}{\scriptstyle 2}
\bigl( 1 \,-\, \frac{\scriptstyle 1}{\scriptstyle q} \bigr) }}
\end{equation}
\mbox{} \vspace{-0.220cm} \\
for each
$ \;\! t > \tilde{t}_0 $,
$ q > 2 $.
This gives,
letting
$ \;\! q \rightarrow \infty $, \\
\mbox{} \vspace{-0.700cm} \\
\begin{equation}
\notag
\|\, \tilde{\mbox{\boldmath $v$}}(\cdot,t)
\,-\, \mbox{\boldmath $v$}(\cdot,t) \,
\|_{\mbox{}_{\scriptstyle L^{\infty}(\mathbb{K})}}
\;\!\leq\:
\frac{\:\Gamma}{\sqrt{\;\!2\;}\;} \,
(\;\! \tilde{t}_0 \!\;\!-\;\! t_0 )^{\mbox{}^{\scriptstyle
\!\:\! \frac{\scriptstyle 1}{\scriptstyle 2} }}
\:\!
\|\, \mbox{\boldmath $u$}_0 \,
\|_{\mbox{}_{\scriptstyle L^{2}(\mathbb{R}^{3})}}
  ^{\mbox{}^{\scriptstyle \:\! 2}}
(\;\! t - \tilde{t}_0 )^{\mbox{}^{\scriptstyle \!\!\!
-\, \frac{\scriptstyle 3}{\scriptstyle 2} }}
\end{equation}
\mbox{} \vspace{-0.220cm} \\
for each
$ \;\! t > \tilde{t}_0 $,
with
$ \mathbb{K} \subset \mathbb{R}^{3} \!\:\!$
compact
{\em arbitrary}.
This estimate
clearly implies (2.16).
}
\mbox{} \hfill $\Box$ \\
%
%
\mbox{} \vspace{-0.650cm} \\

For the next fundamental result
reviewed in this section,
given in Theorem~2.4,
we will need
the following
elementary
Sobolev-Nirenberg-Gagliardo
(\mbox{\small SNG})
inequalities
for
arbitrary
$ {\displaystyle
\;\!
\mbox{u} \in H^{2}(\mathbb{R}^{3})
} $: \\
\mbox{} \vspace{-0.600cm} \\
\begin{equation}
\tag{2.17$a$}
\|\: \mbox{u} \:\|_{\mbox{}_{\scriptstyle \infty}}
\;\!\leq\,
K_{\mbox{}_{\!\;\!0}} \,
\|\: \mbox{u} \:
\|_{\mbox{}_{\scriptstyle L^{2}(\mathbb{R}^{3})}}^{\:\!1/4}
\;\!
\|\: D^{2} \mbox{u} \:
\|_{\mbox{}_{\scriptstyle L^{2}(\mathbb{R}^{3})}}^{\:\!3/4}
\!\:\!,
\qquad
K_{\mbox{}_{\!\;\!0}} \!\;\!<\;\! 0.678,
\end{equation}
\mbox{} \vspace{-0.150cm} \\
see e.g.$\;$\cite{Taylor2011}, Proposition 2.4, p.$\;$13,
or \cite{Schutz2008}, Theorem 4.5.1, p.$\;$52;
$\;\!$and \\
\mbox{} \vspace{-0.600cm} \\
\begin{equation}
\tag{2.17$b$}
\|\, D \:\!\mbox{u} \:
\|_{\mbox{}_{\scriptstyle L^{2}(\mathbb{R}^{3})}}
\leq\,
K_{\mbox{}_{\!\;\!1}} \;\!
\|\: \mbox{u} \:
\|_{\mbox{}_{\scriptstyle L^{2}(\mathbb{R}^{3})}}^{\:\!1/2}
\;\!
\|\: D^{2} \mbox{u} \:
\|_{\mbox{}_{\scriptstyle L^{2}(\mathbb{R}^{3})}}^{\:\!1/2}
\!\!\;\!,
\qquad
K_{\mbox{}_{\!\;\!1}} \!\;\!=\;\! 1,
\end{equation}
\mbox{} \vspace{-0.150cm} \\
which is easily derived
with the Fourier transform.
By (2.17$a$), (2.17$b$),
we then have \\
\mbox{} \vspace{-0.600cm} \\
\begin{equation}
\tag{2.18}
\|\: \mbox{u} \:\|_{\mbox{}_{\scriptstyle \infty}}
\;\!
\|\, D \,\!\mbox{u} \:
\|_{\mbox{}_{\scriptstyle L^{2}(\mathbb{R}^{3})}}^{\:\!1/2}
\leq\,
K_{\mbox{}_{\!\;\!2}} \;\!
\|\: \mbox{u} \:
\|_{\mbox{}_{\scriptstyle L^{2}(\mathbb{R}^{3})}}^{\:\!1/2}
\;\!
\|\: D^{2} \mbox{u} \:
\|_{\mbox{}_{\scriptstyle L^{2}(\mathbb{R}^{3})}}
\!\:\!,
\quad \;\,
K_{\mbox{}_{\!\;\!2}} \!\;\!=\;\!
K_{\mbox{}_{\!\;\!0}} \;\!
K_{\mbox{}_{\!\;\!1}}^{\:\!1/2}
\!\:\!< 1.
\end{equation}
\mbox{} \vspace{-0.700cm} \\
%
%
%
%
%
{\bf Theorem 2.4.}
\textit{%
Let
$ \,\mbox{\boldmath $u$}(\cdot,t) $,
$ t > 0 $,
be any particular
Leray-Hopf's solution to $\;\!(1.1)$.
Then,
there exists
$ \,t_{\!\;\!\ast\ast} \!\;\!\gg\!\;\! 1 $
$(\:\!t_{\!\;\!\ast\ast}\!\:\!$
depending on the solution
$\;\!\mbox{\boldmath $u$}$\/$)$
sufficiently large that
$ {\displaystyle
\|\, D \mbox{\boldmath $u$}(\cdot,t) \,
\|_{\mbox{}_{\scriptstyle L^{2}(\mathbb{R}^{3})}}
\!\:\!
} $
is a smooth, monotonically decreasing
function of $ \, t $
on $\:\! [\,t_{\!\;\!\ast\ast}\!\;\!, \:\!\mbox{\small $\infty$}\:\![\:\! $.
} \\
%
%
%
\mbox{} \vspace{-0.100cm} \\
{\small
{\bf Proof:}
The following argument
is adapted from
\cite{KreissHagstromLorenzZingano2002},
Lemma 2.2.
Let
$ {\displaystyle
\;\!
t_{0}
\!\;\!\geq
t_{\ast}
} $
(to be chosen shortly),
with
$ {\displaystyle
\;\!
t_{\ast}
\!\:\!\gg 1
\;\!
} $
given in (1.3).
Let
$ \;\! t > t_{0} $.
Applying
$ {\displaystyle
D_{\mbox{}_{\scriptstyle \!\ell}} \!\;\!=\:\!
\partial/\partial \:\!x_{\mbox{}_{\scriptstyle \!\ell}}
} $
to the first equation in (1.1$a$),
taking the inner product with
$ {\displaystyle
\:\!D_{\mbox{}_{\scriptstyle \!\ell}}
\mbox{\boldmath $u$}(\cdot,t)
} $
and
integrating on
$ {\displaystyle
\mathbb{R}^{3} \!\times\!\;\!
[\, t_{0}, \;\!t\;\!]
} $,
we get,
summing over $ 1 \leq \ell \leq 3 $, \\
\mbox{} \vspace{-0.200cm} \\
\mbox{} \hspace{+3.500cm}
$ {\displaystyle
\|\, D \mbox{\boldmath $u$}(\cdot,t) \,
\|_{\mbox{}_{\scriptstyle L^{2}(\mathbb{R}^{3})}}^{\:\!2}
+\:
2 \!\!\;\!
\int_{\mbox{\footnotesize $ t_{0} $}}
    ^{\mbox{\footnotesize $\:\!t$}}
\!
\|\, D^{2} \mbox{\boldmath $u$}(\cdot,s) \,
\|_{\mbox{}_{\scriptstyle L^{2}(\mathbb{R}^{3})}}^{\:\!2}
ds
\;\;\!=
} $ \\
\mbox{} \vspace{+0.020cm} \\
\mbox{} \hspace{+1.000cm}
$ {\displaystyle
=\;\:\!
\|\, D \mbox{\boldmath $u$}(\cdot,t_{0}) \,
\|_{\mbox{}_{\scriptstyle L^{2}(\mathbb{R}^{3})}}^{\:\!2}
\!\;\!+\:
2
\sum_{i, \, j, \, \ell}
\int_{\mbox{\footnotesize $ t_{0} $}}
    ^{\mbox{\footnotesize $\:\!t$}}
\!
\int_{\mathbb{R}^{3}}
\!\!\!\;\!
u_{i}(x,s) \;\!
D_{\mbox{}_{\scriptstyle \!\!\;\!\ell}} u_{j}(x,s)
\;\!
D_{\scriptstyle \!j}
D_{\mbox{}_{\scriptstyle \!\!\;\!\ell}}
u_{i}(x,s) \,
dx \: ds
} $ \\
\mbox{} \vspace{+0.050cm} \\
\mbox{} \hspace{+1.000cm}
$ {\displaystyle
\leq\;
\|\, D \mbox{\boldmath $u$}(\cdot,t_{0}) \,
\|_{\mbox{}_{\scriptstyle L^{2}(\mathbb{R}^{3})}}^{\:\!2}
+\:
2 \!\!\:\!
\int_{\mbox{\footnotesize $ t_{0} $}}
    ^{\mbox{\footnotesize $\:\!t$}}
\!
\|\, \mbox{\boldmath $u$}(\cdot,s) \,
\|_{\mbox{}_{\scriptstyle \infty}}
\:\!
\|\, D \mbox{\boldmath $u$}(\cdot,s) \,
\|_{\mbox{}_{\scriptstyle L^{2}(\mathbb{R}^{3})}}
\|\, D^{2} \mbox{\boldmath $u$}(\cdot,s) \,
\|_{\mbox{}_{\scriptstyle L^{2}(\mathbb{R}^{3})}}
\;\!
ds
} $ \\
\mbox{} \vspace{+0.025cm} \\
\mbox{} \hspace{+1.000cm}
$ {\displaystyle
\leq\;
\|\, D \mbox{\boldmath $u$}(\cdot,t_{0}) \,
\|_{\mbox{}_{\scriptstyle L^{2}(\mathbb{R}^{3})}}^{\:\!2}
+\:
2 \!\!\:\!
\int_{\mbox{\footnotesize $ t_{0} $}}
    ^{\mbox{\footnotesize $\:\!t$}}
\!
\|\, \mbox{\boldmath $u$}(\cdot,s) \,
\|_{\mbox{}_{\scriptstyle L^{2}(\mathbb{R}^{3})}}^{\:\!1/2}
\:\!
\|\, D \mbox{\boldmath $u$}(\cdot,s) \,
\|_{\mbox{}_{\scriptstyle L^{2}(\mathbb{R}^{3})}}^{\:\!1/2}
\|\, D^{2} \mbox{\boldmath $u$}(\cdot,s) \,
\|_{\mbox{}_{\scriptstyle L^{2}(\mathbb{R}^{3})}}^{\:\!2}
\;\!
ds
} $, \\
\mbox{} \vspace{+0.100cm} \\
by (2.18),
using
(1.16) and (1.17).
In particular,
we have \\
\mbox{} \vspace{-0.150cm} \\
\mbox{} \hspace{+3.500cm}
$ {\displaystyle
\|\, D \mbox{\boldmath $u$}(\cdot,t) \,
\|_{\mbox{}_{\scriptstyle L^{2}(\mathbb{R}^{3})}}^{\:\!2}
+\:
2 \!\!\;\!
\int_{\mbox{\footnotesize $ t_{0} $}}
    ^{\mbox{\footnotesize $\:\!t$}}
\!
\|\, D^{2} \mbox{\boldmath $u$}(\cdot,s) \,
\|_{\mbox{}_{\scriptstyle L^{2}(\mathbb{R}^{3})}}^{\:\!2}
ds
\;\;\!\leq
} $ \\
\mbox{} \vspace{-0.700cm} \\
\mbox{} \hfill (2.19) \\
\mbox{} \vspace{-0.400cm} \\
\mbox{} \hspace{+0.300cm}
$ {\displaystyle
\leq\;
\|\, D \mbox{\boldmath $u$}
(\cdot,t_{0}) \,
\|_{\mbox{}_{\scriptstyle L^{2}(\mathbb{R}^{3})}}^{\:\!2}
+\:
2 \!
\int_{\mbox{\footnotesize $ t_{0} $}}
    ^{\mbox{\footnotesize $\:\!t$}}
\Bigl[\;
\|\, \mbox{\boldmath $u$}_{0} \,
\|_{\mbox{}_{\scriptstyle L^{2}(\mathbb{R}^{3})}}
\:\!
\|\, D \mbox{\boldmath $u$}(\cdot,s) \,
\|_{\mbox{}_{\scriptstyle L^{2}(\mathbb{R}^{3})}}
\;\!
\Bigr]^{\!1/2}
\|\, D^{2} \mbox{\boldmath $u$}(\cdot,s) \,
\|_{\mbox{}_{\scriptstyle L^{2}(\mathbb{R}^{3})}}^{\:\!2}
ds
} $ \\
\mbox{} \vspace{+0.100cm} \\
for all
$ {\displaystyle
\;\!
t \geq t_{0}
} $.
$\!$We then choose
$ \;\! t_{0} \geq t_{\!\;\!\ast} \;\!$
such that,
by (1.2):
$ {\displaystyle
\;\!
\|\, \mbox{\boldmath $u$}_{0} \,
\|_{\mbox{}_{\scriptstyle L^{2}(\mathbb{R}^{3})}}
\,\!
\|\, D \mbox{\boldmath $u$}(\cdot,t_{0}) \,
\|_{\mbox{}_{\scriptstyle L^{2}(\mathbb{R}^{3})}}
\!\!\;\!< 1
} $. \linebreak
\mbox{} \vspace{-0.550cm} \\
In fact,
with this choice,
it follows from (2.19)
that \\
\mbox{} \vspace{-0.600cm} \\
\begin{equation}
\tag{2.20}
\|\, \mbox{\boldmath $u$}_{0} \,
\|_{\mbox{}_{\scriptstyle L^{2}(\mathbb{R}^{3})}}
\,\!
\|\, D \mbox{\boldmath $u$}(\cdot,s) \,
\|_{\mbox{}_{\scriptstyle L^{2}(\mathbb{R}^{3})}}
<\;\! 1
\qquad
\forall \;\;\!s \geq t_0.
\end{equation}
\mbox{} \vspace{-0.200cm} \\
\mbox{$[\,$}Proof
of (2.20): if false,
there would be
$\;\! t_{1} \!\;\!> t_{0} \;\!$
such that
$ {\displaystyle
\;\!
\|\, \mbox{\boldmath $u$}_{0} \,
\|_{\mbox{}_{\scriptstyle L^{2}(\mathbb{R}^{3})}}
\,\!
\|\, D \mbox{\boldmath $u$}(\cdot,s) \,
\|_{\mbox{}_{\scriptstyle L^{2}(\mathbb{R}^{3})}}
\!< 1
} $ \linebreak
\mbox{} \vspace{-0.530cm} \\
for all
$ \;\!t_0 \leq s < t_{1}$,
while
$ {\displaystyle
\;\!
\|\, \mbox{\boldmath $u$}_{0} \,
\|_{\mbox{}_{\scriptstyle L^{2}(\mathbb{R}^{3})}}
\,\!
\|\, D \mbox{\boldmath $u$}(\cdot,t_{1}) \,
\|_{\mbox{}_{\scriptstyle L^{2}(\mathbb{R}^{3})}}
\!= 1
} $.
Taking $\;\! t = t_{1} $
in (2.19) above, \linebreak
\mbox{} \vspace{-0.530cm} \\
this would give
$ {\displaystyle
\|\, D \mbox{\boldmath $u$}(\cdot,t_1) \,
\|_{\mbox{}_{\scriptstyle L^{2}(\mathbb{R}^{3})}}
\!\!\;\!\leq\!\;\!
\|\, D \mbox{\boldmath $u$}(\cdot,t_0) \,
\|_{\mbox{}_{\scriptstyle L^{2}(\mathbb{R}^{3})}}
\!\:\!
} $,
so that
$ {\displaystyle
\|\, \mbox{\boldmath $u$}_0 \,
\|_{\mbox{}_{\scriptstyle L^{2}(\mathbb{R}^{3})}}
\!\;\!
\|\, D \mbox{\boldmath $u$}(\cdot,t_1) \,
\|_{\mbox{}_{\scriptstyle L^{2}(\mathbb{R}^{3})}}
} $ \\
\mbox{} \vspace{-0.530cm} \\
$ {\displaystyle
\leq\:\!
\|\, \mbox{\boldmath $u$}_0 \,
\|_{\mbox{}_{\scriptstyle L^{2}(\mathbb{R}^{3})}}
\|\, D \mbox{\boldmath $u$}(\cdot,t_0) \,
\|_{\mbox{}_{\scriptstyle L^{2}(\mathbb{R}^{3})}}
\!< 1
} $.
$\!$This contradiction shows (2.20),
as claimed.
$\!$\mbox{\scriptsize QED (2.20)}$]$ \\
%
\mbox{} \vspace{-0.700cm} \\

From (2.19) and (2.20),
it then follows
that \\
\mbox{} \vspace{-0.650cm} \\
\begin{equation}
\tag{2.21}
\|\, D \mbox{\boldmath $u$}(\cdot,t) \,
\|_{\mbox{}_{\scriptstyle L^{2}(\mathbb{R}^{3})}}^{\:\!2}
+\:
2 \, \gamma
\!\!\;\!
\int_{\mbox{\footnotesize $ t_{2} $}}
    ^{\mbox{\footnotesize $\:\!t$}}
\!
\|\, D^{2} \mbox{\boldmath $u$}(\cdot,s) \,
\|_{\mbox{}_{\scriptstyle L^{2}(\mathbb{R}^{3})}}^{\:\!2}
ds
\;\;\!\leq\;
\|\, D \mbox{\boldmath $u$}
(\cdot,t_{2}) \,
\|_{\mbox{}_{\scriptstyle L^{2}(\mathbb{R}^{3})}}^{\:\!2}
\end{equation}
\mbox{} \vspace{-0.150cm} \\
for all
$\;\! t \geq t_{2} \geq t_{0} $,
where
$ {\displaystyle
\,
\gamma \:\! := \:
1 \,-\: \|\, \mbox{\boldmath $u$}(\cdot,t_0) \,
\|_{\mbox{}_{\scriptstyle L^{2}(\mathbb{R}^{3})}}^{\:\!1/2}
\|\, D \mbox{\boldmath $u$}(\cdot,t_0) \,
\|_{\mbox{}_{\scriptstyle L^{2}(\mathbb{R}^{3})}}^{\:\!1/2}
\!
} $
is some positive con-\linebreak
stant.
This shows the
result,
where
$ \;\!t_{\!\;\!\ast\ast} \!\:\!=\;\! t_{0} $
with
$\;\!t_{0} \geq t_{\!\;\!\ast} $
as chosen
in (2.20) above.
}
\mbox{} \hfill $\Box$ \\
%
\mbox{} \vspace{-0.650cm} \\
\nl
{\bf Remark 2.1.}
As shown in
\cite{KreissHagstromLorenzZingano2002},
one can readily obtain
from the proof of Theorem~2.4 \linebreak
that
one has
$ {\displaystyle
\,
t_{\ast\ast} \!\;\!<
0.212 \cdot
\|\, \mbox{\boldmath $u$}_{0} \;\!
\|_{\scriptstyle L^{2}(\mathbb{R}^{3})}
  ^{\:\!4}
\!\;\!
} $
always.
A more elaborated argument
developed here
in the \mbox{\small \sc Appendix}
produces the much sharper estimate \\
\mbox{} \vspace{-0.075cm} \\
\mbox{} \hspace{+3.800cm}
\fbox{%
\begin{minipage}{6.150cm}
\nl
\mbox{} \vspace{-0.400cm} \\
$ {\displaystyle
\mbox{} \;\;
t_{\ast\ast}
\,<\;
0.000\,753\,026 \cdot
\|\, \mbox{\boldmath $u$}_{0} \;\!
\|_{\mbox{}_{\scriptstyle L^{2}(\mathbb{R}^{3})}}
  ^{\mbox{}^{\scriptstyle \:\!4}}
\!\:\!
} $, \\
\end{minipage}
}
\mbox{} \vspace{-1.050cm} \\
\mbox{} \hfill (2.22) \\
\nl
\mbox{} \vspace{-0.225cm} \\
giving a practical upper bound
on how much one should wait
in numerical experiments
before we may witness
any loss of regularity
on the part of $ \mbox{\boldmath $u$}(\cdot,t) $,
if this ever happens ---
in any case,
we have
$ {\displaystyle
\;\!
\mbox{\boldmath $u$} \in C^{\infty}(\mathbb{R}^{3}
\!\!\;\!\times\!\:\! [\, t_{\ast\ast}\!\;\!, \infty))
} $.
$\!\;\!$Other estimates for $ t_{\ast\ast} $
have appeared in the literature,
see e.g.$\;$\cite{Galdi2000, Leray1934, LorenzZingano2003, %
RobinsonSadowskiSilva2012};
(2.22) is
the sharpest of its kind. \linebreak
Whether one can really have
$ {\displaystyle
\mbox{\boldmath $u$} \!\not{\!\!\;\!\in}\;
C^{\infty}(\mathbb{R}^{3}
\!\!\;\!\times\!\:\! (\;\! 0, \infty\:\!))
} $
for some Leray-Hopf's
solutions
is not really known and
remains one of the famous fundamental
open questions regarding the Leray-Hopf's
solutions of the Navier-Stokes equations
\cite{Constantin2001, Fefferman2006, Heywood1990}.

\mbox{} \vspace{-0.800cm} \\

Our final basic result to be collected
in this section for convenience of the reader
is the following
fundamental property,
which is a direct consequence
of (1.2), (2.3) and
the monotonicity of
$ {\displaystyle
\;\!
\|\, D \mbox{\boldmath $u$}(\cdot,t) \,
\|_{\mbox{}_{\scriptstyle L^{2}(\mathbb{R}^{3})}}
\!\:\!
} $
for large $ t $,
as given in Theorem~2.4
above. \\
\nl
%
%
%
%
%
{\bf Theorem 2.5.}
\textit{%
Let
$ \,\mbox{\boldmath $u$}(\cdot,t) $,
$ t > 0 $,
be any particular
Leray-Hopf's solution to $\;\!(1.1)$.
Then
} \\
\mbox{} \vspace{-1.100cm} \\
\begin{equation}
\tag{2.23}
\lim_{\;t\,\rightarrow\,\infty}
\,
t^{\:\!1/2}
\,
\|\, D \mbox{\boldmath $u$}(\cdot,t) \,
\|_{\mbox{}_{\scriptstyle L^{2}(\mathbb{R}^{3})}}
=\; 0.
\end{equation}
%
%
\nl
\mbox{} \vspace{-0.475cm} \\
%
%
{\small
{\bf Proof:}
The following argument
is taken from
\cite{KreissHagstromLorenzZingano2002},
Lemma 2.1.
If (2.23) were false,
there would then
exist an increasing sequence
$ \;\! t_{\ell} \;\!\mbox{\scriptsize $\nearrow$}\;\! \infty $
(with
$ \:\! t_{\ell} \geq t_{\!\;\!\ast\ast} $
and
$ \;\! t_{\ell} \geq 2 \;\! t_{\ell - 1} $
for all $\ell $,
say) \linebreak
and some
fixed
$ \;\!\eta > 0 \;\!$
such that \\
\mbox{} \vspace{-0.620cm} \\
\begin{equation}
\notag
\mbox{} \;\;\;
t_{\ell} \:
\|\, D \mbox{\boldmath $u$}(\cdot,t_{\ell}) \,
\|_{\mbox{}_{\scriptstyle L^{2}(\mathbb{R}^{3})}}^{\:\!2}
\:\!\geq\; \eta
\qquad
\forall \;\;\! \ell.
\end{equation}
\mbox{} \vspace{-0.240cm} \\
In particular,
this would give \\
\mbox{} \vspace{-0.600cm} \\
\begin{equation}
\notag
\int_{\mbox{}_{\mbox{\footnotesize $\!\:\!t_{\ell - 1}$}}}
    ^{\mbox{\footnotesize $\:\!t_{\ell}$}}
\!\!\!\!\!
\|\, D \mbox{\boldmath $u$}(\cdot,s) \,
\|_{\mbox{}_{\scriptstyle L^{2}(\mathbb{R}^{3})}}^{\:\!2}
\:\!dt
\:\,\!\geq\:
(\:\!t_{\ell} -\;\! t_{\ell - 1})
\,
\|\, D \mbox{\boldmath $u$}(\cdot,t_{\ell}) \,
\|_{\mbox{}_{\scriptstyle L^{2}(\mathbb{R}^{3})}}^{\:\!2}
\!\:\!\geq\,
\mbox{\footnotesize $ {\displaystyle \frac{1}{2} }$}
\: t_{\ell} \,
\|\, D \mbox{\boldmath $u$}(\cdot,t_{\ell}) \,
\|_{\mbox{}_{\scriptstyle L^{2}(\mathbb{R}^{3})}}^{\:\!2}
\!\:\!\geq\,
\mbox{\footnotesize $ {\displaystyle \frac{1}{2} }$}
\, \eta
\end{equation}
\mbox{} \vspace{-0.100cm} \\
for all $\ell $,
in contradiction with (1.2), (2.3).
This concludes
the proof of (2.23),
as claimed.
}
\mbox{} \hfill $\Box$ \\
%
%
\nl
\mbox{} \vspace{-1.250cm} \\
\newpage
%

%
%
\mbox{} \vspace{-0.950cm} \\

{\bf 2. Proof of the \mbox{\boldmath $L^{2}\!\:\!$}
results (1.4) and (1.8)} \\
\setcounter{section}{3}
\mbox{} \vspace{-0.500cm} \\

Now that the basic properties
of Leray-Hopf's solutions
given above
have been established,
it becomes much easier to obtain
estimates like
(1.4), (1.5), (1.8) or (1.9).
In this section,
we consider (1.4) and (1.8).
Let
then
$ {\displaystyle
\;\!
\mbox{\boldmath $u$}(\cdot,t)
\in
L^{\infty}(\:\![\;\!0, \infty\:\![\:\!, L^{2}_{\sigma}(\mathbb{R}^{3})\:\!)
} $
$ {\displaystyle
\cap \,
L^{2}(\:\![\;\!0, \infty\:\![\:\!, \stackrel{.}{H}\!\!\mbox{}^{1}(\mathbb{R}^{3})\:\!)
} $
be any such solution
to the initial value problem (1.1$a$),~(1.1$b$),
and let $ \:\!t_{\!\;\!\ast} \!\;\!\gg\!\;\! 1 $
be large enough that (1.3) holds.
Taking
$ \:\!t_0 \!\;\!\geq t_{\!\;\!\ast} \!\:\!$
(arbitrary),
we thus have the representation \\
\mbox{} \vspace{-0.675cm} \\
\begin{equation}
\tag{3.1}
\mbox{\boldmath $u$}(\cdot,t)
\;=\;
e^{\Delta (\mbox{\footnotesize $t$} \;\!-\, \mbox{\footnotesize $t_0$})}
\:\!
\mbox{\boldmath $u$}(\cdot,t_0)
\:+\!\;\!
\int_{\mbox{\footnotesize $\!\;\!t_0$}}
    ^{\mbox{\footnotesize $\:\!t$}}
\!\!\;\!
e^{\Delta (\mbox{\footnotesize $t$} \;\!-\, \mbox{\footnotesize $s$})}
\:\!
\mbox{\boldmath $Q$}(\cdot,s)
\: ds,
\qquad
t \geq t_0,
\end{equation}
\mbox{} \vspace{-0.190cm} \\
by Duhamel's principle,
where
$ {\displaystyle
\;\!
\mbox{\boldmath $Q$}(\cdot,s)
\;\!
} $
is defined in (2.4). \\
\nl
%
%
%
%
%
{\bf Theorem 3.1} ({\em Leray's $L^{2}\!$ conjecture\/}).
\textit{%
One has
} \\
\mbox{} \vspace{-0.750cm} \\
\begin{equation}
\tag{3.2}
\lim_{t\,\rightarrow\,\infty}
\,
\|\, \mbox{\boldmath $u$}(\cdot,t) \,
\|_{\mbox{}_{\scriptstyle L^{2}(\mathbb{R}^{3})}}
=\;0.
\end{equation}
%
\nl
\mbox{} \vspace{-0.450cm} \\
{\small
{\bf Proof:}
Given $ \:\!\epsilon > 0 $,
let $ \;\!t_0 \!\;\!\geq t_{\ast} $
be chosen
large enough so that,
by Theorem 2.5, \\
\mbox{} \vspace{-0.550cm} \\
\begin{equation}
\tag{3.3}
\mbox{} \hspace{+0.500cm}
t^{\:\!1/2} \,
\|\, D \mbox{\boldmath $u$}(\cdot,t) \,
\|_{\mbox{}_{\scriptstyle L^{2}(\mathbb{R}^{3})}}
\leq\: \epsilon
\qquad
\forall \;\;\!
t \geq t_0.
\end{equation}
\mbox{} \vspace{-0.175cm} \\
From the representation (3.1)
for $ \mbox{\boldmath $u$}(\cdot,t) $,
this gives \\
\mbox{} \vspace{-0.075cm} \\
\mbox{} \hspace{+0.500cm}
$ {\displaystyle
\|\, \mbox{\boldmath $u$}(\cdot,t) \,
\|_{\mbox{}_{\scriptstyle L^{2}(\mathbb{R}^{3})}}
\;\!\leq\;
\|\: e^{\Delta (t \;\!-\, t_0)} \:\!
\mbox{\boldmath $u$}(\cdot,t_0) \,
\|_{\mbox{}_{\scriptstyle L^{2}(\mathbb{R}^{3})}}
\:\!+
\int_{\mbox{\scriptsize $t_0$}}
    ^{\mbox{\scriptsize $t$}}
\!
\|\: e^{\Delta (t \;\!-\, s)} \:\!
\mbox{\boldmath $Q$}(\cdot,s) \,
\|_{\mbox{}_{\scriptstyle L^{2}(\mathbb{R}^{3})}}
\:\!
ds
} $ \\
\mbox{} \vspace{-0.050cm} \\
\mbox{} \hfill
$ {\displaystyle
\leq\:
\|\: e^{\Delta (t \;\!-\, t_0)} \:\!
\mbox{\boldmath $u$}(\cdot,t_0) \,
\|_{\mbox{}_{\scriptstyle L^{2}(\mathbb{R}^{3})}}
\:\!+\:
K \!\!
\int_{\mbox{\scriptsize $t_0$}}
    ^{\mbox{\scriptsize $t$}}
\!
(t - s)^{-\,3/4} \:
\|\, \mbox{\boldmath $u$}(\cdot,s) \,
\|_{\mbox{}_{\scriptstyle L^{2}(\mathbb{R}^{3})}}
\:\!
\|\, D \mbox{\boldmath $u$}(\cdot,s) \,
\|_{\mbox{}_{\scriptstyle L^{2}(\mathbb{R}^{3})}}
\:\!
ds
} $ \\
\mbox{} \vspace{-0.050cm} \\
\mbox{} \hfill
$ {\displaystyle
\leq\:
\|\: e^{\Delta (t \;\!-\, t_0)} \:\!
\mbox{\boldmath $u$}(\cdot,t_0) \,
\|_{\mbox{}_{\scriptstyle L^{2}(\mathbb{R}^{3})}}
\:\!+\:
K \,
\|\, \mbox{\boldmath $u$}_0 \,
\|_{\mbox{}_{\scriptstyle L^{2}(\mathbb{R}^{3})}}
\!\!\;\!
\int_{\mbox{\scriptsize $t_0$}}
    ^{\mbox{\scriptsize $t$}}
\!\;\!
(t - s)^{-\,3/4} \:
\|\, D \mbox{\boldmath $u$}(\cdot,s) \,
\|_{\mbox{}_{\scriptstyle L^{2}(\mathbb{R}^{3})}}
\:\!
ds
} $ \\
\mbox{} \vspace{-0.050cm} \\
%
\mbox{} \hspace{+1.275cm}
$ {\displaystyle
\leq\:
\|\: e^{\Delta (\mbox{\scriptsize $t$} \;\!-\, t_0)} \:\!
\mbox{\boldmath $u$}(\cdot,t_0) \,
\|_{\mbox{}_{\scriptstyle L^{2}(\mathbb{R}^{3})}}
\:\!+\:
K \,
\|\, \mbox{\boldmath $u$}_0 \,
\|_{\mbox{}_{\scriptstyle L^{2}(\mathbb{R}^{3})}}
\;\!
\epsilon
\!\:\!
\int_{\mbox{\scriptsize $t_0$}}
    ^{\mbox{\scriptsize $t$}}
\!\;\!
(t - s)^{-\,3/4} \, s^{-\,1/2}
\,ds
} $
\mbox{} \hfill [$\,$by (3.3)$\,$] \\
\mbox{} \vspace{+0.150cm} \\
for all $ \:\!t > t_0 $,
where
$ \:\!K \!\:\!=\:\! (\:\!8 \:\!\pi )^{-\,3/4} \!\:\!$,
and where
we have used (1.2), (2.5$a$).
Observing that \\
\mbox{} \vspace{-0.575cm} \\
\begin{equation}
\notag
\int_{\mbox{\scriptsize $t_0$}}
    ^{\mbox{\scriptsize $t$}}
\!\;\!
(t - s)^{-\,3/4} \, s^{-\,1/2}
\,ds
\;\leq\;
6 \: \sqrt[4]{\;\!2\;}
\qquad
\forall
\;\,
t \;\!\geq\;\! t_0 + 1,
\end{equation}
\mbox{} \vspace{-0.075cm} \\
we then have,
for all $ \:\! t \geq\;\! t_0 +\;\! 1 $: \\
\mbox{} \vspace{-0.650cm} \\
\begin{equation}
\notag
\|\, \mbox{\boldmath $u$}(\cdot,t) \,
\|_{\mbox{}_{\scriptstyle L^{2}(\mathbb{R}^{3})}}
\;\!\leq\;
\|\: e^{\Delta (t \;\!-\, t_0)} \:\!
\mbox{\boldmath $u$}(\cdot,t_0) \,
\|_{\mbox{}_{\scriptstyle L^{2}(\mathbb{R}^{3})}}
\:\!+\:
\|\, \mbox{\boldmath $u$}_0 \,
\|_{\mbox{}_{\scriptstyle L^{2}(\mathbb{R}^{3})}}
\epsilon.
\end{equation}
\mbox{} \vspace{-0.150cm} \\
Recalling that,
for the heat semigroup,
we have
$ {\displaystyle
\;\!
\|\: e^{\Delta (t \;\!-\, t_0)} \:\!
\mbox{\boldmath $u$}(\cdot,t_0) \,
\|_{\mbox{}_{\scriptstyle L^{2}(\mathbb{R}^{3})}}
\!\rightarrow \;\!0
\;\!
} $
as $ \;\!t \rightarrow \infty $, \linebreak
\mbox{} \vspace{-0.550cm} \\
it follows that \\
\mbox{} \vspace{-0.900cm} \\
\begin{equation}
\notag
\|\, \mbox{\boldmath $u$}(\cdot,t) \,
\|_{\mbox{}_{\scriptstyle L^{2}(\mathbb{R}^{3})}}
\;\!\leq\;
(\;\! 1 \,+\,
\|\, \mbox{\boldmath $u$}_0 \,
\|_{\mbox{}_{\scriptstyle L^{2}(\mathbb{R}^{3})}})
\: \epsilon
\end{equation}
\mbox{} \vspace{-0.250cm} \\
for all $ \;\! t \gg 1 $.
Since
$ \;\!\epsilon > 0 \;\!$
is arbitrary,
this shows (3.2),
completing our argument.
}
\mbox{} $\Box$ \\
%
\nl
%
%
%
%
{\bf Theorem 3.2.}
\textit{%
Given any $\;\!t_0 \geq 0 $,
one has
} \\
\mbox{} \vspace{-0.700cm} \\
\begin{equation}
\tag{3.4}
\lim_{t\,\rightarrow\,\infty}
\,
t^{\:\!1/4} \;\!
\|\, \mbox{\boldmath $u$}(\cdot,t) \;\!-\,
e^{\Delta (\:\! \mbox{\footnotesize $t$} \,-\, \mbox{\footnotesize $t_0$})}
\:\! \mbox{\boldmath $u$}(\cdot,t_0) \,
\|_{\mbox{}_{\scriptstyle L^{2}(\mathbb{R}^{3})}}
\:\!=\;
0.
\end{equation}
%
%
\nl
\mbox{} \vspace{-0.500cm} \\
%
%
{\small
{\bf Proof:}
By Theorem 2.2,
it is sufficient to show (3.4)
in the case
$ \:\!t_0 \!\;\!\geq t_{\!\;\!\ast} \!\;\!$,
where (3.1) holds.
Given $ \;\!\epsilon > 0 $,
let
$ \:\!t_{\!\;\!\epsilon} \!\:\!> t_{0} \!\;\!$
be large enough that
we have,
using Theorem 2.5
again, \\
\mbox{} \vspace{-0.575cm} \\
\begin{equation}
\tag{3.5}
\mbox{} \hspace{+0.500cm}
t^{\:\!1/2} \,
\|\, D \mbox{\boldmath $u$}(\cdot,t) \,
\|_{\mbox{}_{\scriptstyle L^{2}(\mathbb{R}^{3})}}
\leq\: \epsilon
\qquad
\forall \;\;\!
t \geq t_{\!\;\!\epsilon}.
\end{equation}
\mbox{} \vspace{-0.225cm} \\
By (3.1)
and (1.2), (2.5$a$),
we then get \\
\mbox{} \vspace{-0.125cm} \\
\mbox{} \hspace{+0.500cm}
$ {\displaystyle
t^{\:\!1/4} \,
\|\, \mbox{\boldmath $u$}(\cdot,t) \,-\:
e^{\Delta (t \;\!-\, t_0)} \:\!
\mbox{\boldmath $u$}(\cdot,t_0) \,
\|_{\mbox{}_{\scriptstyle L^{2}(\mathbb{R}^{3})}}
\:\!\leq\;
t^{\:\!1/4}
\!\!\!\:\!
\int_{\mbox{\scriptsize $\!\;\!t_0$}}
    ^{\mbox{\scriptsize $t$}}
\!\!\:\!
\|\: e^{\Delta (t \;\!-\, s)} \:\!
\mbox{\boldmath $Q$}(\cdot,s) \,
\|_{\mbox{}_{\scriptstyle L^{2}(\mathbb{R}^{3})}}
\:\!
ds
} $ \\
\mbox{} \vspace{-0.050cm} \\
\mbox{} \hfill
$ {\displaystyle
\leq\;
I(\:\!t, \:\!t_{\!\;\!\epsilon})
\;+\:
K \;\! t^{\:\!1/4}
\!\!\!\:\!
\int_{\mbox{\scriptsize $\!\;\!t_{\epsilon}$}}
    ^{\mbox{\scriptsize $t$}}
\!
(t - s)^{-\,3/4} \,
\|\, \mbox{\boldmath $u$}(\cdot,s) \,
\|_{\mbox{}_{\scriptstyle L^{2}(\mathbb{R}^{3})}}
\:\!
\|\, D \mbox{\boldmath $u$}(\cdot,s) \,
\|_{\mbox{}_{\scriptstyle L^{2}(\mathbb{R}^{3})}}
\:\!
ds
} $ \\
\mbox{} \vspace{-0.050cm} \\
\mbox{} \hspace{+3.450cm}
$ {\displaystyle
\leq\;
I(\:\!t, \:\!t_{\!\;\!\epsilon})
\;+\:
K \,
\|\, \mbox{\boldmath $u$}_0 \,
\|_{\mbox{}_{\scriptstyle L^{2}(\mathbb{R}^{3})}}
\epsilon
\:
t^{\:\!1/4}
\!\!\!\:\!
\int_{\mbox{\scriptsize $\!\;\!t_{\epsilon}$}}
    ^{\mbox{\scriptsize $t$}}
\!
(t - s)^{-\,3/4} \, s^{-\,1/2}
\,ds
} $
\mbox{} \hfill [$\,$by (3.5)$\,$] \\
\mbox{} \vspace{-0.050cm} \\
\mbox{} \hspace{+3.450cm}
$ {\displaystyle
\leq\;
I(\:\!t, \:\!t_{\!\;\!\epsilon})
\;+\:
0.636 \:
\|\, \mbox{\boldmath $u$}_0 \,
\|_{\mbox{}_{\scriptstyle L^{2}(\mathbb{R}^{3})}}
\:\!
\epsilon
\;
t^{\:\!1/4}
\;\!
(\:\! t - t_{\!\;\!\epsilon})^{-\,1/4}
} $ \\
\mbox{} \vspace{+0.025cm} \\
for all
$ \;\! t > t_{\!\;\!\epsilon}$,
where
$ \:\!K \!\:\!=\:\! (\:\!8 \:\!\pi )^{-\,3/4} \!\:\!$
and \\
\mbox{} \vspace{-0.100cm} \\
\mbox{} \hspace{+2.130cm}
$ {\displaystyle
I(\:\!t, \:\!t_{\!\;\!\epsilon})
\;=\:
K \;\! t^{\:\!1/4}
\!\!\!\:\!
\int_{\mbox{\scriptsize $\!\;\!t_{0}$}}
    ^{\mbox{\scriptsize $t_{\epsilon}$}}
\!
(t - s)^{-\,3/4} \,
\|\, \mbox{\boldmath $u$}(\cdot,s) \,
\|_{\mbox{}_{\scriptstyle L^{2}(\mathbb{R}^{3})}}
\:\!
\|\, D \mbox{\boldmath $u$}(\cdot,s) \,
\|_{\mbox{}_{\scriptstyle L^{2}(\mathbb{R}^{3})}}
\:\!
ds
} $ \\
\mbox{} \vspace{-0.050cm} \\
\mbox{} \hspace{+3.450cm}
$ {\displaystyle
\leq\:
K \,
t^{\:\!1/4}
\;\!
(\:\!t - t_{\!\;\!\epsilon})^{-\,3/4}
\!\!\!\:\!
\int_{\mbox{\scriptsize $\!\;\!t_{0}$}}
    ^{\mbox{\scriptsize $t_{\epsilon}$}}
\!
\|\, \mbox{\boldmath $u$}(\cdot,s) \,
\|_{\mbox{}_{\scriptstyle L^{2}(\mathbb{R}^{3})}}
\:\!
\|\, D \mbox{\boldmath $u$}(\cdot,s) \,
\|_{\mbox{}_{\scriptstyle L^{2}(\mathbb{R}^{3})}}
\:\!
ds
} $. \\
\mbox{} \vspace{-0.120cm} \\
Therefore,
we obtain \\
\mbox{} \vspace{-0.550cm} \\
\begin{equation}
\notag
t^{\:\!1/4} \,
\|\, \mbox{\boldmath $u$}(\cdot,t) \,-\:
e^{\Delta (t \;\!-\, t_0)} \:\!
\mbox{\boldmath $u$}(\cdot,t_0) \,
\|_{\mbox{}_{\scriptstyle L^{2}(\mathbb{R}^{3})}}
\:\!\leq\;
(\;\! 1 \,+\,
\|\, \mbox{\boldmath $u$}_0 \,
\|_{\mbox{}_{\scriptstyle L^{2}(\mathbb{R}^{3})}})
\: \epsilon
\end{equation}
\mbox{} \vspace{-0.300cm} \\
for all $ \;\! t \gg 1 $.
%
%
This gives (3.4),
and our $L^2$ discussion
is now complete,
as claimed.
}
\mbox{} \hfill $\Box$ \\
%
\nl

%
%
\nl
\mbox{} \vspace{-2.000cm} \\

{\bf 4. Proof of the supnorm results (1.5) and (1.9)} \\
\setcounter{section}{4}
\mbox{} \vspace{-0.500cm} \\

In this section,
we follow a similar path
to obtain the more delicate
supnorm estimates (1.5) and (1.9).
Let then
$ {\displaystyle
\;\!
\mbox{\boldmath $u$}(\cdot,t)
\in
L^{\infty}(\:\![\;\!0, \infty\:\![\:\!, L^{2}_{\sigma}(\mathbb{R}^{3})\:\!)
} $
$ {\displaystyle
\cap \,
L^{2}(\:\![\;\!0, \infty\:\![\:\!, \stackrel{.}{H}\!\!\mbox{}^{1}(\mathbb{R}^{3})\:\!)
} $
be any given Leray-Hopf's solution
to the Cauchy problem (1.1).
Again,
we take advantage
of the strong regularity
properties of
$ \mbox{\boldmath $u$}(\cdot,t) $
for $ t \geq t_{\ast} \gg 1 $
(see (1.3) above), \linebreak
using the representation (3.1)
and
the fundamental results
(2.5$b$), (2.23) and (3.2)
already obtained. \\
\nl
%
%
%
%
%
{\bf Theorem 4.1.}
\textit{%
One has
} \\
\mbox{} \vspace{-0.900cm} \\
\begin{equation}
\tag{4.1}
\lim_{t\,\rightarrow\,\infty}
\, t^{\:\!3/4} \,
\|\, \mbox{\boldmath $u$}(\cdot,t) \,
\|_{\mbox{}_{\scriptstyle L^{\infty}(\mathbb{R}^{3})}}
=\;0.
\end{equation}
\nl
%
%
{\small
{\bf Proof:}
Given $ \:\!0 < \epsilon \leq 1/2 $,
let $ \;\!t_0 \!\;\!\geq t_{\ast} $
be large enough that,
by (2.23) and (3.2) above,
we have \\
\mbox{} \vspace{-0.900cm} \\
\begin{equation}
\tag{4.2$a$}
\mbox{} \hspace{+0.500cm}
t^{\:\!1/2} \,
\|\, D \mbox{\boldmath $u$}(\cdot,t) \,
\|_{\mbox{}_{\scriptstyle L^{2}(\mathbb{R}^{3})}}
\:\!\leq\; \epsilon
\qquad
\forall \;\;\!
t \geq t_0
\end{equation}
\mbox{} \vspace{-0.600cm} \\
and \\
\mbox{} \vspace{-1.150cm} \\
\begin{equation}
\tag{4.2$b$}
\mbox{} \hspace{+0.500cm}
\|\, \mbox{\boldmath $u$}(\cdot,t) \,
\|_{\mbox{}_{\scriptstyle L^{2}(\mathbb{R}^{3})}}
\:\!\leq\; \epsilon
\qquad
\forall \;\;\!
t \geq t_0.
\end{equation}
\mbox{} \vspace{-0.100cm} \\
From the representation (3.1)
for $ \mbox{\boldmath $u$}(\cdot,t) $,
we obtain,
by (2.5$b$)
and (4.2$a$), \\
\mbox{} \vspace{-0.100cm} \\
\mbox{} \hspace{+0.500cm}
$ {\displaystyle
\|\, \mbox{\boldmath $u$}(\cdot,t) \,
\|_{\mbox{}_{\scriptstyle \infty}}
\;\!\leq\;
\|\: e^{\Delta (t \;\!-\, t_0)} \:\!
\mbox{\boldmath $u$}(\cdot,t_0) \,
\|_{\mbox{}_{\scriptstyle \infty}}
\:\!+
\int_{\mbox{\scriptsize $\!\;\!t_0$}}
    ^{\mbox{\scriptsize $\:\!t$}}
\!
\|\: e^{\Delta (t \;\!-\, s)} \:\!
\mbox{\boldmath $Q$}(\cdot,s) \,
\|_{\mbox{}_{\scriptstyle \infty}}
\;\!
ds
} $ \\
\mbox{} \vspace{-0.050cm} \\
\mbox{} \hfill
$ {\displaystyle
\leq\:
\|\: e^{\Delta (t \;\!-\, t_0)} \:\!
\mbox{\boldmath $u$}(\cdot,t_0) \,
\|_{\mbox{}_{\scriptstyle \infty}}
\:\!+\:
K \!\!
\int_{\mbox{\scriptsize $\!\;\!t_0$}}
    ^{\mbox{\scriptsize $\:\!t$}}
\!\:\!
(t - s)^{-\,3/4} \,
\|\, \mbox{\boldmath $u$}(\cdot,s) \,
\|_{\mbox{}_{\scriptstyle \infty}}
\;\!
\|\, D \mbox{\boldmath $u$}(\cdot,s) \,
\|_{\mbox{}_{\scriptstyle L^{2}(\mathbb{R}^{3})}}
\;\!
ds
} $ \\
%
%
\mbox{} \vspace{-0.050cm} \\
\mbox{} \hspace{+2.450cm}
$ {\displaystyle
\leq\;
K \:\! (\:\!t - t_0)^{-\,3/4} \;\!
\|\, \mbox{\boldmath $u$}(\cdot,t_0) \,
\|_{\mbox{}_{\scriptstyle L^{2}(\mathbb{R}^{3})}}
\:\!+\;
K \:\! \epsilon
\!\!\;\!
\int_{\mbox{\scriptsize $\!\;\!t_0$}}
    ^{\mbox{\scriptsize $\:\!t$}}
\!\:\!
(t - s)^{-\,3/4} \;\!
s^{-\,1/2} \,
\|\, \mbox{\boldmath $u$}(\cdot,s) \,
\|_{\mbox{}_{\scriptstyle \infty}}
\;\!
ds
} $ \\
%
%
%
%
\mbox{} \vspace{+0.050cm} \\
for all $ \:\!t > t_0 $,
where
$ \:\!K \!\:\!=\:\! (\:\!8 \:\!\pi )^{-\,3/4} \!\:\!$.
That is, \\
\mbox{} \vspace{-0.100cm} \\
\mbox{} \hspace{-0.250cm}
$ {\displaystyle
\|\, \mbox{\boldmath $u$}(\cdot,t) \,
\|_{\mbox{}_{\scriptstyle \infty}}
\:\!\leq\:
K \:\! (\:\!t - t_0)^{-\,3/4} \;\!
\|\, \mbox{\boldmath $u$}(\cdot,t_0) \,
\|_{\mbox{}_{\scriptstyle L^{2}(\mathbb{R}^{3})}}
\!\;\!+\;\!
K \:\! \epsilon
\!\!\;\!
\int_{\mbox{\scriptsize $\!\;\!t_0$}}
    ^{\mbox{\scriptsize $\:\!t$}}
\!\:\!
(t - s)^{-\,3/4} \;\!
s^{-\,1/2} \,
\|\, \mbox{\boldmath $u$}(\cdot,s) \,
\|_{\mbox{}_{\scriptstyle \infty}}
\;\!
ds
} $
\hfill (4.3) \\
\mbox{} \vspace{-0.000cm} \\
for all $\;\!t > t_0 $.
We claim that
this gives \\
\mbox{} \vspace{-0.570cm} \\
\begin{equation}
\tag{4.4}
t^{\:\!3/4} \,
\|\, \mbox{\boldmath $u$}(\cdot,t) \,
\|_{\mbox{}_{\scriptstyle \infty}}
\;\!<\; \epsilon
\qquad
\forall \;\;\!
t \;\!\geq\;\! 2\,(\:\!t_0 + 1),
\end{equation}
\mbox{} \vspace{-0.190cm} \\
which implies (4.1).
In what follows,
we will prove (4.4) above.
Given any
$ {\displaystyle
\;\!
\hat{t} \:\!\geq\;\! 2\,(\:\!t_0 + 1)
} $
(fixed, but otherwise arbitrary),
let
$ \;\! \hat{t}_{1} \!:=\:\! \hat{t}/2 $.
Setting \\
\mbox{} \vspace{-0.600cm} \\
\begin{equation}
\tag{4.5}
{\sf U}(t) :=\:
(\:\!t - \hat{t}_{1}\!\;\!)^{3/4} \,
\|\, \mbox{\boldmath $u$}(\cdot,t) \,
\|_{\mbox{}_{\scriptstyle \infty}}
\!\;\!,
\qquad
t \:\!\geq\;\! \hat{t}_{1},
\end{equation}
%
%
we obtain,
applying (4.3)
[$\;\!$with $ t_0 = \hat{t}_1 $ there$\;\!$]$\:\!$: \\
\mbox{} \vspace{-0.100cm} \\
\mbox{} \hspace{+0.740cm}
$ {\displaystyle
{\sf U}(t)
\;\leq\;
K \,
\|\, \mbox{\boldmath $u$}(\cdot,\hat{t}_1) \,
\|_{\mbox{}_{\scriptstyle L^{2}(\mathbb{R}^{3})}}
\:\!+\;
K \:\! \epsilon \:
(\:\!t - \hat{t}_{1}\!\;\!)^{3/4}
\!\!
\int_{\mbox{\scriptsize $\hat{t}_1$}}
    ^{\mbox{\scriptsize $\:\!t$}}
\!\:\!
(\:\!t - s)^{-\,3/4} \;\!
s^{-\,1/2} \,
\|\, \mbox{\boldmath $u$}(\cdot,s) \,
\|_{\mbox{}_{\scriptstyle \infty}}
\;\!
ds
} $ \\
\mbox{} \vspace{-0.050cm} \\
\mbox{} \hspace{+1.690cm}
$ {\displaystyle
=\;
K \,
\|\, \mbox{\boldmath $u$}(\cdot,\hat{t}_1) \,
\|_{\mbox{}_{\scriptstyle L^{2}(\mathbb{R}^{3})}}
\:\!+\;
K \:\! \epsilon \:
(\:\!t - \hat{t}_{1}\!\;\!)^{3/4}
\!\!
\int_{\mbox{\scriptsize $\hat{t}_1$}}
    ^{\mbox{\scriptsize $\:\!t$}}
\!\:\!
(\:\!t - s)^{-\,3/4} \;\!
s^{-\,1/2} \;\!
(s - \hat{t}_{1}\!\;\!)^{-\,3/4} \:
{\sf U}(s) \,
ds
} $ \\
\mbox{} \vspace{-0.050cm} \\
\mbox{} \hfill
$ {\displaystyle
\leq\;
K \,
\|\, \mbox{\boldmath $u$}(\cdot,\hat{t}_1) \,
\|_{\mbox{}_{\scriptstyle L^{2}(\mathbb{R}^{3})}}
\:\!+\;
K \:\! \epsilon \;\:\!
\hat{t}_{1}^{\,-\,1/2} \;\!
(\:\!t - \hat{t}_{1}\!\;\!)^{3/4}
\!\!
\int_{\mbox{\scriptsize $\hat{t}_1$}}
    ^{\mbox{\scriptsize $\:\!t$}}
\!\:\!
(\:\!t - s)^{-\,3/4} \,
(s - \hat{t}_{1}\!\;\!)^{-\,3/4} \:
{\sf U}(s) \,
ds
} $ \\
\mbox{} \vspace{+0.075cm} \\
for $\;\! t \geq \hat{t}_{1} $,
so that,
by (4.2$b$),
we have \\
\mbox{} \vspace{-0.000cm} \\
\mbox{} \hspace{+0.750cm}
$ {\displaystyle
{\sf U}(t)
\;\leq\;
K \:\! \epsilon
\;+\:
K \:\! \epsilon \;\:\!
\hat{t}_{1}^{\,-\,1/2} \;\!
(\:\!t - \hat{t}_{1}\!\;\!)^{3/4}
\!\!
\int_{\mbox{\scriptsize $\hat{t}_1$}}
    ^{\mbox{\scriptsize $\:\!t$}}
\!\:\!
(\:\!t - s)^{-\,3/4} \,
(s - \hat{t}_{1}\!\;\!)^{-\,3/4} \:
{\sf U}(s) \,
ds
} $,
\mbox{} \hspace{+0.450cm}
$ t \;\!\geq\;\! \hat{t}_{1}$, \\
\mbox{} \vspace{+0.075cm} \\
where
$ \:\!K \!\:\!=\:\! (\:\!8 \:\!\pi )^{-\,3/4} \!\;\!$.
$\!$Thus,
setting
$ {\displaystyle
\;\!
\hat{\mathbb{U}}
\!\;\!:=\,
\max\:\{\, {\sf U}(t) \!\:\!:\: \hat{t}_{1} \leq \:\!t \leq \:\! \hat{t} \,\}
} $,
we get,
for
$\;\! \hat{t}_{1} \leq \:\! t \leq \:\! \hat{t} $: \linebreak
\mbox{} \vspace{-0.000cm} \\
\mbox{} \hspace{+1.850cm}
$ {\displaystyle
{\sf U}(t)
\;\leq\;
K \:\! \epsilon
\;+\:
K \:\! \epsilon \;\:\!
\hat{t}_{1}^{\,-\,1/2} \;\!
(\:\!t - \hat{t}_{1}\!\;\!)^{3/4} \;
\hat{\mathbb{U}}
\!
\int_{\mbox{\scriptsize $\hat{t}_1$}}
    ^{\mbox{\scriptsize $\:\!t$}}
\!\:\!
(\:\!t - s)^{-\,3/4} \,
(s - \hat{t}_{1}\!\;\!)^{-\,3/4} \;\!
ds
} $ \\
\mbox{} \vspace{-0.000cm} \\
\mbox{} \hspace{+2.750cm}
$ {\displaystyle
\leq\;
K \:\! \epsilon
\;+\;
8 \, \sqrt{\;\!2\;} \,
K \:\! \epsilon \;\:\!
\hat{t}_{1}^{\,-\,1/2} \;\!
(\:\!t - \hat{t}_{1}\!\;\!)^{1/4} \;
\hat{\mathbb{U}}
} $ \\
\mbox{} \vspace{-0.000cm} \\
\mbox{} \hspace{+2.750cm}
$ {\displaystyle
\leq\;
K \:\! \epsilon
\;+\;
8 \, \sqrt{\;\!2\;} \,
K \:\! \epsilon \;\:\!
\hat{t}_{1}^{\,-\,1/4} \;
\hat{\mathbb{U}}
} $. \\
\mbox{} \vspace{-0.000cm} \\
Recalling that
$ \;\!\hat{t}_{1} \geq 1 $,
$ \epsilon \;\!\leq 1/2 $,
we then get
$ {\displaystyle
\,
{\sf U}(t) \;\!\leq\;\!
K \epsilon
\;\!+\;\!
8 \, \sqrt{\;\!2\;} \,
K \epsilon \;
\hat{\mathbb{U}}
\,\leq\;\!
K \epsilon
\;\!+\;\!
0.504 \; \hat{\mathbb{U}}
} $,
for each
$ \;\! \hat{t}_{1} \!\;\!\leq \;\! t \:\!\leq\;\! \hat{t} $.
This gives
$ {\displaystyle
\,
\hat{\mathbb{U}}
\,\leq\;\!
K \epsilon
\;\!+\;\!
0.504 \; \hat{\mathbb{U}}
} $,
that is,
$ {\displaystyle
\;\!
\hat{\mathbb{U}}
\:\!<\;\!
0.180 \: \epsilon
} $.
In particular, \\
\mbox{} \vspace{-0.600cm} \\
\begin{equation}
\notag
\Bigl(\;\! \frac{\:\!\hat{t}\:\!}{\mbox{\footnotesize $2$}}
\;\!\Bigr)^{\!3/4} \,
\|\, \mbox{\boldmath $u$}(\cdot,\hat{t}\:\!)\,
\|_{\mbox{}_{\scriptstyle \infty}}
\;\!=\:
{\sf U}(\:\!\hat{t}\:\!)
\,\leq\,
\hat{\mathbb{U}}
\,<\:
0.180 \: \epsilon,
\end{equation}
\mbox{} \vspace{-0.200cm} \\
so that
$ {\displaystyle
\;\!
\hat{t}^{\,3/4} \;\!
\|\, \mbox{\boldmath $u$}(\cdot,\hat{t}\:\!)\,
\|_{\mbox{}_{\scriptstyle \infty}}
\!< \epsilon
} $,
where
$ \;\!\hat{t} \;\!\geq\;\! 2 \,(\:\!t_0 + 1 ) \:\!$
is arbitrary.
This shows (4.4),
as claimed,
and the proof of (4.1)
is now complete.
}
\mbox{} \hfill $\Box$ \linebreak
%
\mbox{} \vspace{-0.625cm} \\

In a similar way,
(1.9) can be obtained,
as shown next. \\
\nl
%
%
%
%
{\bf Theorem 4.2.}
\textit{%
Given any $\;\!t_0 \geq 0 $,
one has
} \\
\mbox{} \vspace{-0.750cm} \\
\begin{equation}
\tag{4.6}
\lim_{t\,\rightarrow\,\infty}
\,
t \:
\|\, \mbox{\boldmath $u$}(\cdot,t) \;\!-\,
e^{\Delta (\:\! \mbox{\footnotesize $t$} \,-\, \mbox{\footnotesize $t_0$})}
\:\! \mbox{\boldmath $u$}(\cdot,t_0) \,
\|_{\mbox{}_{\scriptstyle L^{\infty}(\mathbb{R}^{3})}}
\:\!=\;
0.
\end{equation}
%
%
\nl
\mbox{} \vspace{-0.500cm} \\
%
%
{\small
{\bf Proof:}
By Theorem 2.2,
it is sufficient to show (4.6)
in the case
$ \:\!t_0 \!\;\!\geq t_{\!\;\!\ast} \!\;\!$,
where (3.1) holds.
Given $ \;\!0 < \epsilon \:\!\leq\:\! 1 $,
let
$ \:\!t_{\!\;\!\epsilon} \!\:\!> t_{0} $
be large enough that
we have,
from Theorems 2.5
and 4.1, \\
\mbox{} \vspace{-0.625cm} \\
\begin{equation}
\tag{4.7$a$}
\mbox{} \hspace{+0.500cm}
t^{\:\!1/2} \,
\|\, D \mbox{\boldmath $u$}(\cdot,t) \,
\|_{\mbox{}_{\scriptstyle L^{2}(\mathbb{R}^{3})}}
\;\!\leq\; \epsilon
\qquad
\forall \;\;\!
t \geq t_{\!\;\!\epsilon},
\end{equation}
%
%
and

\mbox{} \vspace{-1.150cm} \\
\begin{equation}
\tag{4.7$b$}
\mbox{} \hspace{+0.500cm}
t^{\:\!3/4} \,
\|\, \mbox{\boldmath $u$}(\cdot,t) \,
\|_{\mbox{}_{\scriptstyle \infty}}
\;\!\leq\; \epsilon
\qquad
\forall \;\;\!
t \geq t_{\!\;\!\epsilon}.
\end{equation}

\mbox{} \vspace{-0.300cm} \\
From (3.1),
we have \\
\mbox{} \vspace{-0.150cm} \\
\mbox{} \hspace{+2.000cm}
$ {\displaystyle
t \;
\|\, \mbox{\boldmath $u$}(\cdot,t) \,-\:
e^{\Delta (t \;\!-\, t_0)} \:\!
\mbox{\boldmath $u$}(\cdot,t_0) \,
\|_{\mbox{}_{\scriptstyle \infty}}
\;\!\leq\;
t
\!\!\;\!
\int_{\mbox{\scriptsize $\!\;\!t_0$}}
    ^{\mbox{\scriptsize $t$}}
\!\!\:\!
\|\: e^{\Delta (t \;\!-\, s)} \:\!
\mbox{\boldmath $Q$}(\cdot,s) \,
\|_{\mbox{}_{\scriptstyle \infty}}
\;\!
ds
} $ \\
\mbox{} \vspace{-0.050cm} \\
\mbox{} \hspace{+7.900cm}
$ {\displaystyle
\leq\;
J_{1}(t) \:+\: J_{2}(t) \:+\: J_{3}(t)
} $
\mbox{} \hfill (4.8) \\
\mbox{} \vspace{-0.300cm} \\
for all
$ \;\! t > t_{\!\;\!\epsilon} $,
where \\
\mbox{} \vspace{-0.150cm} \\
\mbox{} \hspace{+1.450cm}
$ {\displaystyle
J_{1}(t) \;=\;\;\!
t
\!
\int_{\mbox{\scriptsize $\!\;\!t_0$}}
    ^{\mbox{\scriptsize $t_{\epsilon}$}}
\!\!\:\!
\|\: e^{\Delta (t \;\!-\, s)} \:\!
\mbox{\boldmath $Q$}(\cdot,s) \,
\|_{\mbox{}_{\scriptstyle \infty}}
\;\!
ds
} $ \\
\mbox{} \vspace{-0.050cm} \\
\mbox{} \hspace{+2.500cm}
$ {\displaystyle
\leq\;\:\!
(\:\!4 \:\!\pi )^{-\,3/4}
\;
t
\!
\int_{\mbox{\scriptsize $\!\;\!t_0$}}
    ^{\mbox{\scriptsize $t_{\epsilon}$}}
\!\!\:\!
(\:\! t - s )^{-\,3/4} \:
\|\: e^{\;\!\frac{\scriptscriptstyle 1}{\scriptscriptstyle 2}
        \:\!\Delta (t \;\!-\, s)} \:\!
\mbox{\boldmath $Q$}(\cdot,s) \,
\|_{\mbox{}_{\scriptstyle L^{2}(\mathbb{R}^{3})}}
\:\!
ds
} $ \\
\mbox{} \vspace{-0.050cm} \\
\mbox{} \hspace{+2.500cm}
$ {\displaystyle
\leq\;\:\!
(\:\!4 \:\!\pi )^{-\,3/2}
\;
t \!
\int_{\mbox{\scriptsize $\!\;\!t_0$}}
    ^{\mbox{\scriptsize $t_{\epsilon}$}}
\!\!\:\!
(\:\! t - s )^{-\,3/2}
\:
\|\, \mbox{\boldmath $u$}(\cdot,s) \,
\|_{\mbox{}_{\scriptstyle L^{2}(\mathbb{R}^{3})}}
\|\, D \mbox{\boldmath $u$}(\cdot,s) \,
\|_{\mbox{}_{\scriptstyle L^{2}(\mathbb{R}^{3})}}
\:\!
ds
} $ \\
\mbox{} \vspace{-0.050cm} \\
\mbox{} \hspace{+2.500cm}
$ {\displaystyle
\leq\;\:\!
(\:\!4 \:\!\pi )^{-\,3/2}
\;
t \;\:\!
(\:\! t - t_{\!\;\!\epsilon} )^{-\,3/2}
\!\!
\int_{\mbox{\scriptsize $\!\;\!t_0$}}
    ^{\mbox{\scriptsize $t_{\epsilon}$}}
\!\!\:\!
\|\, \mbox{\boldmath $u$}(\cdot,s) \,
\|_{\mbox{}_{\scriptstyle L^{2}(\mathbb{R}^{3})}}
\|\, D \mbox{\boldmath $u$}(\cdot,s) \,
\|_{\mbox{}_{\scriptstyle L^{2}(\mathbb{R}^{3})}}
\:\!
ds
} $,
\mbox{} \hfill (4.9$a$) \\
%
%
\mbox{} \vspace{+0.000cm} \\
by (2.5$a$),
and \\
\mbox{} \vspace{-0.150cm} \\
\mbox{} \hspace{+1.450cm}
$ {\displaystyle
J_{2}(t) \;=\;\;\!
t
\!
\int_{\mbox{\scriptsize $\!\;\! t_{\epsilon}$}}
    ^{\mbox{\scriptsize $\:\!\mu(t)$}}
\!\!\!\!\!\!
\|\: e^{\Delta (t \;\!-\, s)} \:\!
\mbox{\boldmath $Q$}(\cdot,s) \,
\|_{\mbox{}_{\scriptstyle \infty}}
\;\!
ds
} $ \\
\mbox{} \vspace{-0.050cm} \\
\mbox{} \hspace{+2.500cm}
$ {\displaystyle
\leq\;\:\!
(\:\!4 \:\!\pi )^{-\,3/4}
\;
t
\!
\int_{\mbox{\scriptsize $\!\;\! t_{\epsilon}$}}
    ^{\mbox{\scriptsize $\:\!\mu(t)$}}
\!\!\!\!\!\!\;\!
(\:\! t - s )^{-\,3/4} \:
\|\: e^{\;\!\frac{\scriptscriptstyle 1}{\scriptscriptstyle 2}
        \:\!\Delta (t \;\!-\, s)} \:\!
\mbox{\boldmath $Q$}(\cdot,s) \,
\|_{\mbox{}_{\scriptstyle L^{2}(\mathbb{R}^{3})}}
\:\!
ds
} $ \\
\mbox{} \vspace{-0.050cm} \\
\mbox{} \hspace{+2.500cm}
$ {\displaystyle
\leq\;\:\!
(\:\!4 \:\!\pi )^{-\,3/2}
\;
t \!
\int_{\mbox{\scriptsize $\!\;\! t_{\epsilon}$}}
    ^{\mbox{\scriptsize $\:\!\mu(t)$}}
\!\!\!\!\!\!\;\!
(\:\! t - s )^{-\,3/2}
\:
\|\, \mbox{\boldmath $u$}(\cdot,s) \,
\|_{\mbox{}_{\scriptstyle L^{2}(\mathbb{R}^{3})}}
\|\, D \mbox{\boldmath $u$}(\cdot,s) \,
\|_{\mbox{}_{\scriptstyle L^{2}(\mathbb{R}^{3})}}
\:\!
ds
} $ \\
\mbox{} \vspace{-0.050cm} \\
\mbox{} \hspace{+2.500cm}
$ {\displaystyle
\leq\;\:\!
(\:\!2 \:\!\pi )^{-\,3/2}
\;
t \;\:\!
(\:\! t - t_{\!\;\!\epsilon} )^{-\,3/2}
\, \epsilon^{2}
\!\!
\int_{\mbox{\scriptsize $\!\;\! t_{\epsilon}$}}
    ^{\mbox{\scriptsize $\:\!\mu(t)$}}
\!\!\!\!\!\!\;\!
s^{-\,1/2} \;\! ds
\;\;\!<\;\;\!
0.090 \;\:\! \epsilon \;
t \: (\:\! t - t_{\!\;\!\epsilon})^{-\,1}
\!
} $,
\mbox{} \hfill (4.9$b$) \\
\mbox{} \vspace{+0.100cm} \\
by (2.5$a$) and (4.7),
where
$ \;\! \mu(t) \:\!=\;\! (\:\!t + t_{\!\;\!\epsilon})/2 $, \\
\mbox{} \vspace{-0.100cm} \\
\mbox{} \hspace{+1.450cm}
$ {\displaystyle
J_{3}(t) \;=\;\;\!
t
\!
\int_{\mbox{\scriptsize $\!\;\! \mu(t)$}}
    ^{\mbox{\scriptsize $\:\!t$}}
\!\!\!
\|\: e^{\Delta (t \;\!-\, s)} \:\!
\mbox{\boldmath $Q$}(\cdot,s) \,
\|_{\mbox{}_{\scriptstyle \infty}}
\;\!
ds
} $ \\
\mbox{} \vspace{-0.050cm} \\
\mbox{} \hspace{+2.500cm}
$ {\displaystyle
\leq\;\:\!
(\:\!8 \:\!\pi )^{-\,3/4}
\;
t \!
\int_{\mbox{\scriptsize $\!\;\! \mu(t)$}}
    ^{\mbox{\scriptsize $\:\!t$}}
\!\!\!\:\!
(\:\! t - s )^{-\,3/4}
\:
\|\, \mbox{\boldmath $u$}(\cdot,s) \,
\|_{\mbox{}_{\scriptstyle \infty}}
\:\!
\|\, D \mbox{\boldmath $u$}(\cdot,s) \,
\|_{\mbox{}_{\scriptstyle L^{2}(\mathbb{R}^{3})}}
\:\!
ds
} $ \\
\mbox{} \vspace{-0.050cm} \\
\mbox{} \hspace{+2.500cm}
$ {\displaystyle
\leq\;\:\!
(\:\!8 \:\!\pi )^{-\,3/4}
\;
t \; \epsilon^{2}
\!\!\!
\int_{\mbox{\scriptsize $\!\;\! \mu(t)$}}
    ^{\mbox{\scriptsize $\:\!t$}}
\!\!\!\:\!
(\:\! t - s )^{-\,3/4}
\;\!
s^{-\,5/4}
\;\!ds
\;<\;\:\!
0.713 \;\:\!
\epsilon \;
t \;
(\:\! t - t_{\!\;\!\epsilon} )^{-\,1}
\!
} $,
\mbox{} \hfill (4.9$c$) \\
\mbox{} \vspace{+0.000cm} \\
by (2.5$b$) and (4.7).
Therefore,
from (4.8) and (4.9)
above,
we obtain \\
\mbox{} \vspace{-0.575cm} \\
\begin{equation}
\tag{4.10}
t \;
\|\, \mbox{\boldmath $u$}(\cdot,t) \,-\:
e^{\Delta (t \;\!-\, t_0)} \:\!
\mbox{\boldmath $u$}(\cdot,t_0) \,
\|_{\mbox{}_{\scriptstyle \infty}}
\;\!<\;
\epsilon
\qquad
\forall \;\;\!
t \gg 1.
\end{equation}
\mbox{} \vspace{-0.275cm} \\
This shows (4.6),
since
$ {\displaystyle
\;\!
\epsilon \in\;
]\,0, 1 \;\!]
\;\!
} $
is arbitrary.
}
\mbox{} \hfill $\Box$

%

%
\mbox{} \vspace{-1.400cm} \\
%
%
%
%

\mbox{} \hspace{+5.500cm}
\mbox{\large {\sc Appendix}} \\

{\small
Here we show how to obtain
the estimate (2.22)
given in Section 2.
The starting point is
the following inequality,  \\
\mbox{} \vspace{-0.700cm} \\
\begin{equation}
\tag{A.1}
\int_{\mbox{}_{\scriptstyle \mathbb{R}^{3}}}
\!\!\;\!\Bigl\{\!\!\!\!\:\!
\sum_{\mbox{}\;\;\,i, \,j, \,\ell \,=\,1}^{3}
\!\!\!\!
|\, D_{\ell} \;\!u_{i} \;\!| \:
|\, D_{\ell} \;\!u_{j} \;\!| \:
|\, D_{j} \:\!u_{i} \;\!|
\,\Bigr\}
\;\!dx
\:\leq\:
K_{\mbox{}_{\!\:\!3}}^{3} \;\!
\|\, D \:\!\mbox{\boldmath $u$} \,
\|_{\mbox{}_{\scriptstyle L^{2}(\mathbb{R}^{3})}}
  ^{3/2}
\;\!
\|\, D^{2} \mbox{\boldmath $u$} \,
\|_{\mbox{}_{\scriptstyle L^{2}(\mathbb{R}^{3})}}
  ^{3/2}
\!\:\!,
\end{equation}
\mbox{} \vspace{-0.100cm} \\
where
$ \:\!K_{\mbox{}_{\!\:\!3}} \!< 0.581\,862\,001\,307 $
(see \cite{Agueh2008}, Theorem 2.1)
is the constant
in the
Gagliardo-Niren\-berg inequality
$ {\displaystyle
\;\!
\|\: \mbox{u} \:
\|_{\mbox{}_{\scriptstyle L^{3}(\mathbb{R}^{3})}}
\!\leq
K_{\mbox{}_{\!\:\!3}}
\;\!
\|\: \mbox{u} \:
\|_{\mbox{}_{\scriptstyle L^{2}(\mathbb{R}^{3})}}^{1/2}
\|\, D \:\!\mbox{u} \,
\|_{\mbox{}_{\scriptstyle L^{2}(\mathbb{R}^{3})}}^{1/2}
\!\:\!
} $.
$\!$\mbox{[}$\;\!$In fact,
by repeated application
of the Cauchy-Schwarz inequality,
we get \\
\mbox{} \vspace{-0.850cm} \\
\begin{equation}
\notag
\int_{\mbox{}_{\scriptstyle \mathbb{R}^{3}}}
\!\!\;\!\Bigl\{\!\!\!\!\:\!
\sum_{\mbox{}\;\;\,i, \,j, \,\ell \,=\,1}^{3}
\!\!\!\!
|\, D_{\ell} \;\!u_{i} \;\!| \:
|\, D_{\ell} \;\!u_{j} \;\!| \:
|\, D_{j} \:\!u_{i} \;\!|
\,\Bigr\}
\;\!dx
\;\leq\:
\|\: \mbox{v} \:
\|_{\mbox{}_{\scriptstyle L^{3}(\mathbb{R}^{3})}}^{3}
\!,
\quad \;\,
\mbox{v}(x) :=
\Bigl\{\!\!\;\!
\sum_{\mbox{} \;i, \,j, \,=\,1}^{3}
\!\!
|\, D_{j} \:\!u_{i} \;\!|^{\:\!2}
\,
\Bigr\}^{\!\!\;\!1/2}
\!\!\!\!.
\end{equation}
\mbox{} \vspace{-0.100cm} \\
This gives \\
\mbox{} \vspace{-0.425cm} \\
\mbox{} \hspace{+1.500cm}
$ {\displaystyle
\int_{\mbox{}_{\scriptstyle \mathbb{R}^{3}}}
\!\!\;\!\Bigl\{\!\!\!\!\:\!
\sum_{\mbox{}\;\;\,i, \,j, \,\ell \,=\,1}^{3}
\!\!\!\!
|\, D_{\ell} \;\!u_{i} \;\!| \:
|\, D_{\ell} \;\!u_{j} \;\!| \:
|\, D_{j} \:\!u_{i} \;\!|
\,\Bigr\}
\;\!dx
\;\leq\:
K_{\mbox{}_{\!\:\!3}}^{3}
\;\!
\|\: \mbox{v} \:
\|_{\mbox{}_{\scriptstyle L^{2}(\mathbb{R}^{3})}}^{3/2}
\|\, D \:\!\mbox{v} \,
\|_{\mbox{}_{\scriptstyle L^{2}(\mathbb{R}^{3})}}^{3/2}
} $ \\
\mbox{} \vspace{-0.150cm} \\
\mbox{} \hspace{+7.865cm}
$ {\displaystyle
\leq\;\!
K_{\mbox{}_{\!\:\!3}}^{3}
\;\!
\|\, D \:\!\mbox{\boldmath $u$} \,
\|_{\mbox{}_{\scriptstyle L^{2}(\mathbb{R}^{3})}}^{3/2}
\|\, D^{2} \mbox{\boldmath $u$} \,
\|_{\mbox{}_{\scriptstyle L^{2}(\mathbb{R}^{3})}}^{3/2}
} $ \\
\mbox{} \vspace{-0.100cm} \\
by (1.16), as claimed.$\;\!$\mbox{]}
Now,
consider
$ \hat{t} > 0 $
satisfying \\
\mbox{} \vspace{-0.650cm} \\
\begin{equation}
\tag{A.2}
\hat{t} \;>\;
\frac{\;\!1\;\!}{2} \:
K_{\mbox{}_{\!\:\!3}}^{12}
\,
\|\, \mbox{\boldmath $u$}_{0} \;\!
\|_{\mbox{}_{\scriptstyle L^{2}(\mathbb{R}^{3})}}^{4}
\!\:\!:
\end{equation}
\mbox{} \vspace{-0.250cm} \\
Because
(by (1.2))
$ {\displaystyle
\!
\int_{0}^{\;\!\mbox{\footnotesize $\hat{t}$}}
\!\!\;\!
\|\, D \:\!\mbox{\boldmath $u$}(\cdot,t) \,
\|_{\mbox{}_{\scriptstyle L^{2}(\mathbb{R}^{3})}}^{2}
\!\;\!\leq\;\!
\mbox{\footnotesize $ {\displaystyle \frac{1}{2} }$}
\,
\|\, \mbox{\boldmath $u$}_0 \;\!
\|_{\mbox{}_{\scriptstyle L^{2}(\mathbb{R}^{3})}}^{2}
\!
} $, 
there exists
$ \;\!t^{\prime} \!\in (\:\!0, \;\!\hat{t}\,] \:\!$
so that \\
\mbox{} \vspace{-0.500cm} \\
\begin{equation}
\tag{A.3}
\|\, D \:\!\mbox{\boldmath $u$}(\cdot,t^{\prime}) \,
\|_{\mbox{}_{\scriptstyle L^{2}(\mathbb{R}^{3})}}
\leq\,
\|\, \mbox{\boldmath $u$}_{0} \;\!
\|_{\mbox{}_{\scriptstyle L^{2}(\mathbb{R}^{3})}}
\!\cdot\;\!
\frac{\mbox{} \;1\,}{\!\sqrt{\:\!2 \:\! \hat{t}^{\mbox{}}\:}\;}
\!.
\end{equation}
\mbox{} \vspace{-0.140cm} \\
Hence,
by (A.2),
we have
$ {\displaystyle
\:\!
K_{\mbox{}_{\!\:\!3}}^{3}
\;\!
\|\, \mbox{\boldmath $u$}(\cdot,s) \,
\|_{\mbox{}_{\scriptstyle L^{2}(\mathbb{R}^{3})}}^{1/2}
\|\, D\:\!\mbox{\boldmath $u$}(\cdot,s) \,
\|_{\mbox{}_{\scriptstyle L^{2}(\mathbb{R}^{3})}}^{1/2}
\!\!\;\!< 1
\;\!
} $
for all $ \;\!s \geq t^{\prime} $
close to $ t^{\prime} \!\:\!$. \linebreak
\mbox{} \vspace{-0.520cm} \\
From (1.1),
we then get \\
\mbox{} \vspace{-0.150cm} \\
\mbox{} \hspace{+3.500cm}
$ {\displaystyle
\|\, D \mbox{\boldmath $u$}(\cdot,t) \,
\|_{\mbox{}_{\scriptstyle L^{2}(\mathbb{R}^{3})}}^{\:\!2}
+\:
2 \!\!\;\!
\int_{\mbox{\footnotesize $ t^{\prime} $}}
    ^{\mbox{\footnotesize $\:\!t$}}
\!
\|\, D^{2} \mbox{\boldmath $u$}(\cdot,s) \,
\|_{\mbox{}_{\scriptstyle L^{2}(\mathbb{R}^{3})}}^{\:\!2}
ds
} $ \\
\mbox{} \vspace{+0.020cm} \\
\mbox{} \hspace{+1.000cm}
$ {\displaystyle
\leq\;\:\!
\|\, D \mbox{\boldmath $u$}(\cdot,t^{\prime}) \,
\|_{\mbox{}_{\scriptstyle L^{2}(\mathbb{R}^{3})}}^{\:\!2}
\!\;\!+\:
2
\sum_{i, \, j, \, \ell}
\int_{\mbox{\footnotesize $ t^{\prime} $}}
    ^{\mbox{\footnotesize $\:\!t$}}
\!
\int_{\mathbb{R}^{3}}
\!\!\!\;\!
|\, D_{\ell} \:\! u_{i}(x,s) \,|
\;
|\, D_{\ell} \:\! u_{j}(x,s) \,|
\;
|\, D_{\scriptstyle \!j} \:\!u_{i}(x,s) \,|
\;
dx \: ds
} $ \\
\mbox{} \vspace{+0.050cm} \\
\mbox{} \hspace{+1.000cm}
$ {\displaystyle
\leq\;
\|\, D \mbox{\boldmath $u$}(\cdot,t^{\prime}) \,
\|_{\mbox{}_{\scriptstyle L^{2}(\mathbb{R}^{3})}}^{\:\!2}
+\:
2 \!\!\;\!
\int_{\mbox{\footnotesize $ t^{\prime} $}}
    ^{\mbox{\footnotesize $\:\!t$}}
\!
K_{\mbox{}_{\!\:\!3}}^{3}
\:\!
\|\, D \mbox{\boldmath $u$}(\cdot,s) \,
\|_{\mbox{}_{\scriptstyle L^{2}(\mathbb{R}^{3})}}^{\:\!3/2}
\|\, D^{2} \mbox{\boldmath $u$}(\cdot,s) \,
\|_{\mbox{}_{\scriptstyle L^{2}(\mathbb{R}^{3})}}^{\:\!3/2}
\,\!
ds
} $ \\
\mbox{} \vspace{+0.025cm} \\
\mbox{} \hspace{-0.100cm}
$ {\displaystyle
\leq\;
\|\, D \mbox{\boldmath $u$}(\cdot,t^{\prime}) \,
\|_{\mbox{}_{\scriptstyle L^{2}(\mathbb{R}^{3})}}^{\:\!2}
+\:
2 \!\!\;\!
\int_{\mbox{\footnotesize $ t^{\prime} $}}
    ^{\mbox{\footnotesize $\:\!t$}}
\!
\bigl[\,
K_{\mbox{}_{\!\:\!3}}^{3}
\,
\|\, \mbox{\boldmath $u$}(\cdot,s) \,
\|_{\mbox{}_{\scriptstyle L^{2}(\mathbb{R}^{3})}}^{\:\!1/2}
\|\, D \:\!\mbox{\boldmath $u$}(\cdot,s) \,
\|_{\mbox{}_{\scriptstyle L^{2}(\mathbb{R}^{3})}}^{\:\!1/2}
\:\!
\bigr]
\;
\|\, D^{2} \mbox{\boldmath $u$}(\cdot,s) \,
\|_{\mbox{}_{\scriptstyle L^{2}(\mathbb{R}^{3})}}^{2}
ds
} $ \\
\mbox{} \vspace{-0.300cm} \\
\mbox{} \hfill (A.4) \\
\mbox{} \vspace{-0.450cm} \\
for all $ \;\!t \geq t^{\prime} $
close to $ t^{\prime} \!\:\!$.
As in the proof of Theorem 2.4,
this gives \\
\mbox{} \vspace{-0.650cm} \\
\begin{equation}
\tag{A.5}
K_{\mbox{}_{\!\:\!3}}^{3}
\,
\|\, \mbox{\boldmath $u$}(\cdot,t) \,
\|_{\mbox{}_{\scriptstyle L^{2}(\mathbb{R}^{3})}}^{1/2}
\:\!
\|\, D\:\!\mbox{\boldmath $u$}(\cdot,t) \,
\|_{\mbox{}_{\scriptstyle L^{2}(\mathbb{R}^{3})}}^{1/2}
<\, 1,
\qquad
\forall \;\,
t \geq t^{\prime}
\end{equation}
\mbox{} \vspace{-0.200cm} \\
and,
in particular,
as in (A.4) above,
we have \\
\mbox{} \vspace{-0.100cm} \\
\mbox{} \hspace{+3.500cm}
$ {\displaystyle
\|\, D \mbox{\boldmath $u$}(\cdot,t) \,
\|_{\mbox{}_{\scriptstyle L^{2}(\mathbb{R}^{3})}}^{\:\!2}
+\:
2 \!\!\;\!
\int_{\mbox{\footnotesize $ t_0 $}}
    ^{\mbox{\footnotesize $\:\!t$}}
\!
\|\, D^{2} \mbox{\boldmath $u$}(\cdot,s) \,
\|_{\mbox{}_{\scriptstyle L^{2}(\mathbb{R}^{3})}}^{\:\!2}
ds
} $ \\
\mbox{} \vspace{+0.025cm} \\
\mbox{} \hfill
$ {\displaystyle
\leq\;
\|\, D \mbox{\boldmath $u$}(\cdot,t_0) \,
\|_{\mbox{}_{\scriptstyle L^{2}(\mathbb{R}^{3})}}^{\:\!2}
\!\;\!+\:
2 \!\!\;\!
\int_{\mbox{\footnotesize $ t_0 $}}
    ^{\mbox{\footnotesize $\:\!t$}}
\!
\bigl[\,
K_{\mbox{}_{\!\:\!3}}^{3}
\,
\|\, \mbox{\boldmath $u$}(\cdot,s) \,
\|_{\mbox{}_{\scriptstyle L^{2}(\mathbb{R}^{3})}}^{\:\!1/2}
\|\, D \:\!\mbox{\boldmath $u$}(\cdot,s) \,
\|_{\mbox{}_{\scriptstyle L^{2}(\mathbb{R}^{3})}}^{\:\!1/2}
\:\!
\bigr]
\;
\|\, D^{2} \mbox{\boldmath $u$}(\cdot,s) \,
\|_{\mbox{}_{\scriptstyle L^{2}(\mathbb{R}^{3})}}^{2}
ds
} $ \\
\mbox{} \vspace{+0.075cm} \\
for any
$ \;\! t > t_0 \geq t^{\prime} \!\;\!$,
that is,
$ {\displaystyle
\|\, D\:\! \mbox{\boldmath $u$}(\cdot,t) \,
\|_{\mbox{}_{\scriptstyle L^{2}(\mathbb{R}^{3})}}
\!
} $
is monotonically decreasing
in
$ {\displaystyle
[\,t^{\prime}\!, \:\!\infty\:\!)
\supseteq
[\,\hat{t}, \:\!\infty\:\!)
} $,
so that,
by Leray's theory,
$ \mbox{\boldmath $u$} $
must be $C^{\infty}\!$
for $ t > t^{\prime} \!\:\!$.
Recalling (A.2),
this completes
the proof of (2.22),
since
$ {\displaystyle
\:\!
1/2 \cdot
K_{\mbox{}_{\!\:\!3}}^{12}
\!\;\!<\:\!
0.000\,753\,026
} $.
Summarizing it all,
we have shown: \\
\nl
%
%
%
%
{\bf Theorem A.1.}
\textit{%
Let $\;\!\mbox{\boldmath $u$}_0 \!\;\!\in L^{2}_{\sigma}(\mathbb{R}^{3}) $,
and let
$ \;\!\mbox{\boldmath $u$}(\cdot,t) $
be any Leray-Hopf's solution of
$\:\!(1.1)$. \linebreak
Then
there exists
$\,0 \leq t_{\ast\ast} \!\:\!< 0.000\,753\,026 \:
\|\, \mbox{\boldmath $u$}_0 \;\! \|_{L^{2}(\mathbb{R}^{3})}^{\:\!4} \!$
such that
$ \;\!\mbox{\boldmath $u$} \in C^{\infty}(\:\!\mathbb{R}^{3} \!\!\;\!\times\!\:\!
[\,t_{\ast\ast} \!\;\!, \infty)) $
and
$ {\displaystyle
\;\!
\|\, D\:\!\mbox{\boldmath $u$}(\cdot,t) \,
\|_{\mbox{}_{\scriptstyle L^{2}(\mathbb{R}^{3})}}
\!
} $
is finite and monotonically decreasing
everywhere in
$ [\,t_{\ast\ast}\!\;\!, \infty) $.
} \\

%
%
}


%
%

%
%

\nl
\mbox{} \vspace{-0.250cm} \\
\nl
{\small
\begin{minipage}[t]{10.00cm}
\mbox{\normalsize \textsc{Lineia Sch\"utz}} \\
Departamento de Matem\'atica Pura e Aplicada \\
Universidade Federal do Rio Grande do Sul \\
Porto Alegre, RS 91509-900, Brazil \\
E-mail: {\sf lineia.schutz@ufrgs.br}, \\
\mbox{} \hspace{+1.150cm}
        {\sf lineiaschutz@yahoo.com.br} \\
\end{minipage}
\nl
\mbox{} \vspace{-0.450cm} \\
\nl
\begin{minipage}[t]{10.00cm}
\mbox{\normalsize \textsc{Jana\'\i na Pires Zingano}} \\
Departamento de Matem\'atica Pura e Aplicada \\
Universidade Federal do Rio Grande do Sul \\
Porto Alegre, RS 91509-900, Brazil \\
E-mail: {\sf janaina.zingano@ufrgs.br}, \\
\mbox{} \hspace{+1.150cm}
        {\sf jzingano@gmail.com} \\
\end{minipage}
\nl
\mbox{} \vspace{-0.450cm} \\
\nl
\begin{minipage}[t]{10.00cm}
\mbox{\normalsize \textsc{Paulo Ricardo de Avila Zingano}} \\
Departamento de Matem\'atica Pura e Aplicada \\
Universidade Federal do Rio Grande do Sul \\
Porto Alegre, RS 91509-900, Brazil \\
E-mail: {\sf paulo.zingano@ufrgs.br}, \\
\mbox{} \hspace{+1.150cm}
        {\sf zingano@gmail.com} \\
\end{minipage}
}
%
%

%
%

\end{document}